\documentclass[12pt,a4paper,reqno]{article}%
\usepackage{amssymb}
\usepackage{amsmath}
\usepackage{amsfonts}
\usepackage{graphicx}
\usepackage{rotating}
\usepackage{mathrsfs}
\usepackage[all]{xy}
\usepackage{hyperref}
\SelectTips{cm}{}
\newtheorem{thm}{Theorem}
\newtheorem{lem}[thm]{Lemma}
\newtheorem{rem}{Remark}
\newtheorem{defn}{Definition}
\newtheorem{exmp}{Example}
\newtheorem{prop}[thm]{Proposition}
\newtheorem{cor}[thm]{Corollary}
\newtheorem{assum}{Assumption}
\newenvironment{pf}[1][Proof]{\noindent\textbf{#1.} }{\ \rule{0.5em}{0.5em}}
\title{\bf The Lagrangian-Hamiltonian Formalism for Higher Order Field Theories}
\author{\sc{L.~Vitagliano}\thanks{{\bf e}-{\it mail}: \texttt{lvitagliano@unisa.it}}\\
\small{DMI, Universit\`a degli Studi di Salerno, and}\\ \small{Istituto Nazionale di Fisica Nucleare, GC Salerno}\\
\small{Via Ponte don Melillo, 84084 Fisciano (SA), Italy}}

\begin{document}
\maketitle
\begin{abstract}
We generalize the Lagrangian-Hamiltonian formalism of Skinner and Rusk to higher order field theories on fiber bundles. As a byproduct we solve the long standing problem of defining, in a coordinate free manner, a Hamiltonian formalism for higher order Lagrangian field theories. Namely, our formalism does only depend on the action functional and, therefore, unlike previously proposed ones, is free from any relevant ambiguity.
\end{abstract}

\newpage

\section*{Introduction}

First order Lagrangian mechanics can be generalized to higher order Lagrangian
field theory. Moreover, the latter has got a very elegant geometric (and
homological) formulation (see, for instance, \cite{v84}) on which there is
general consensus. On the other hand, it seems that the generalization of
Hamiltonian mechanics of Lagrangian systems to higher order field theory
presents some more problems. Several answers have been proposed (see, for
instance, \cite{d35,d77,aa82,s82,k84b,sc90,s91,k02,av04} and the references
therein) to the question: is there any reasonable, higher order, field
theoretic analogue of Hamiltonian mechanics? In our opinion, none of them is
satisfactorily natural, especially because of the common emergence of
ambiguities due to either the arbitrary choice of a coordinate system
\cite{d35} or the choice of a Legendre transform \cite{sc90,s91,av04}. Namely,
the latter seems not to be uniquely definable, except in the case of first
order Lagrangian field theories when a satisfactory Hamiltonian formulation
can be presented in terms of multisymplectic geometry (see, for instance,
\cite{gs73} - see also \cite{r05} for a recent review, and the references therein).

Nevertheless, it is still desirable to have a Hamiltonian formulation of
higher order Lagrangian field theories enjoying the same nice properties as
Hamiltonian mechanics, which 1) is natural, i.e., is independent of the choice
of any structure other than the action functional, 2) gives rise to first
order equations of motion, 3) takes advantage of the (pre-)symplectic
geometry of the phase space, 4) is a natural starting point for gauge reduction,
5) is a natural starting point for quantization. The relationship between the Euler-Lagrange equations and the Hamilton
equations deserves a special mention. The Legendre transform maps injectively solutions of the former to
solutions of the latter, but, generically, Hamilton equations are not
equivalent to Euler-Lagrange ones \cite{gs73}. However, the difference between
the two is a pure gauge and, therefore, it is irrelevant from a physical point
of view.

In this paper we achieve the goal of finding a natural (in the above mentioned
sense), geometric, higher order, field theoretic analogue of Hamiltonian
mechanics of Lagrangian systems in two steps: first, we find a higher order,
field theoretic analogue of the Skinner and Rusk \textquotedblleft mixed
Lagrangian-Hamiltonian\textquotedblright\ formalism \cite{s83,sr83,sr83b} (see also \cite{e...04}), which is rather straightforward (see \cite{c...09} for a different, finite dimensional approach, to the same problem) and, second, we show
that the derived theory \textquotedblleft projects to a smaller
space\textquotedblright\ which is naturally interpreted as phase space. Local
expressions of the field equations on the phase space are nothing but de
Donder equations \cite{d35} and, therefore, are naturally interpreted as the
higher order, field theoretic, coordinate free analogue of Hamilton equations.
A central role is played in the paper by multisymplectic geometry in the form
of partial differential (PD, in the following) Hamiltonian system theory,
which has been developed in \cite{v09}.

The paper is divided into nine sections. The first four sections contain
reviews of the main aspects of the geometry underlying the paper. They have
been included in order to make the paper as self-consistent as possible. The
next five sections contain most of the original results.

The first section summarizes the notations and conventions adopted throughout the
paper. It also contains references to some differential geometric facts which
are often used in the subsequent sections. Finally, in Section \ref{SecConvSkinRusk} we briefly review the Skinner-Rusk formalism \cite{sr83}. Section \ref{SecGeomPDE} is a short review of the geometric theory of partial
differential equations (PDEs) (see, for instance, \cite{b...99}). Section
\ref{SecCartLagr} outlines the properties of the main geometric structure of jet
spaces and PDEs, the Cartan distribution, and reviews the geometric
formulation of the calculus of variations \cite{v84}. Section \ref{SecPDHam}
reviews the theory of PD-Hamiltonian systems and their PD-Hamilton equations
\cite{v09}. Moreover, it contains examples of morphisms of PDEs coming from
such theory. These examples are presented here for the first time.

In Section \ref{SecLH} we present the higher order, field theoretic analogue
of Skinner-Rusk mixed Lagrangian-Hamiltonian formalism for mechanics. In
Section \ref{SecLH} we also discuss the relationship between the field equations
in the Lagrangian-Hamiltonian formalism (now on, ELH equations) and the
Euler-Lagrange equations. In Section \ref{SecTransfELH} we discuss some
natural transformations of the ELH equations. As a byproduct, we prove that
they are independent of the choice of a Lagrangian density, in the class of
those yielding the same Euler-Lagrange equations, up to isomorphisms. ELH
equations are, therefore, as natural as possible. In Section \ref{SecHamForm}
we present our proposal for a Hamiltonian, higher order, field theory. Since
we don't use any additional structure other than the ELH equations and the
order of a Lagrangian density, we judge our theory satisfactorily natural.
Moreover, the associated field equations (HDW equations) are first order and,
more specifically, of the PD-Hamilton kind. In Section \ref{SecLegTransf} we
study the relationship between the HDW equations and the Euler-Lagrange
equations. As a byproduct, we derive a new (and, in our opinion,
satisfactorily natural) definition of Legendre transform for higher order,
Lagrangian field theories. It is a non-local morphism of the Euler-Lagrange
equations into the HDW equations. Finally, in Section \ref{SecKdV} we apply the theory to the KdV equation which can be derived from a second order variational principle.

\section{Notations, Conventions and the Skinner-Rusk
Formalism\label{SecConvSkinRusk}}

In this section we collect notations and conventions about some general
constructions in differential geometry that will be used in the following.

Let $N$ be a smooth manifold. If $L\subset N$ is a submanifold, we denote by
$i_{L}:L\hookrightarrow N$ the inclusion. We denote by $C^{\infty}(N)$ the
$\mathbb{R}$--algebra of smooth, $\mathbb{R}$--valued functions on $N$. We
will always understand a vector field $X$ on $N$ as a derivation $X:C^{\infty
}(N)\longrightarrow C^{\infty}(N)$. We denote by $\mathrm{D}(N)$ the
$C^{\infty}(N)$--module of vector fields over $N$, by $\Lambda(M)=\bigoplus
_{k}\Lambda^{k}(N)$ the graded $\mathbb{R}$--algebra of differential forms
over $N$ and by $d:\Lambda(N)\longrightarrow\Lambda(N)$ the de Rham
differential. If $F:N_{1}\longrightarrow N$ is a smooth map of manifolds, we
denote by $F^{\ast}:\Lambda(N)\longrightarrow\Lambda(N_{1})$ the pull--back via $F$.
We will understand everywhere the wedge product $\wedge$ of differential
forms, i.e., for $\omega,\omega_{1}\in\Lambda(N)$, we will write $\omega\omega_{1}$ instead of
$\omega\wedge\omega_{1}$, .

Let $\alpha:A\longrightarrow N$ be an affine bundle (for instance, a vector
bundle) and $F:N_{1}\longrightarrow N$ a smooth map of manifolds. Let
$\mathscr{A}$ be the affine space of smooth sections of $\alpha$. For
$a\in\mathscr{A}$ and $x\in N$ we put, sometimes, $a_{x}:=a(x)$. The affine
bundle on $N_{1}$ induced by $\alpha$ via $F$ will be denoted by $F^{\circ
}(\alpha):F^{\circ}(A)\longrightarrow N$:
\[%
\begin{array}
[c]{c}%
\xymatrix{F^\circ(A) \ar[r] \ar[d]_-{F^\circ(\alpha)} & A \ar[d]^-{\alpha} \\ N_1 \ar[r]^-F & N }
\end{array}
,
\]
and the space of its sections by $F^{\circ}(\mathscr{A})$. For any section
$a\in\mathscr{A}$ there exists a unique section, which we denote by $F^{\circ
}(a)\in F^{\circ}(\mathscr{A})$, such that the diagram
\[
\xymatrix{F^\circ(A) \ar[r]  & A  \\
N_1 \ar[r]^-F    \ar[u]^-{F^\circ(a)}                    &  N \ar[u]_-{a}}
\]
commutes. If $F:N_{1}\longrightarrow N$ is the embedding of a submanifold, we
also write $\bullet\ |_{F}$ (or, simply, $\bullet\ |_{N_{1}}$) for $F^{\circ
}({}\bullet{})$, and refer to it as the restriction of \textquotedblleft%
${}\bullet{}$\textquotedblright\ to $N_{1}$ (via $F$), whatever the object
\textquotedblleft${}\bullet{}$\textquotedblright\ is (an affine bundle, its
total space, its space of sections or a section of it).

We will always understand the sum over repeated upper-lower (multi)indexes. Our notations about
multiindexes are the following. We will use the capital letters $I,J,K$ for multiindexes. Let $n$ be a positive integer. A multiindex of length $k$ is a $k$tuple of indexes $I=(i_1,\ldots,i_k)$, $i_1,\ldots,i_k \leq n$. We identify multiindexes differing only by the order of the entries. If $I$ is a multiindex of lenght $k$, we put $|I| := k$. Let $I=(i_1,\ldots,i_k)$ and $J=(j_1,\ldots,j_l)$ be multiindexes, and $i$ an index. We denote by $IJ$ (resp.{} $Ii$) the multiindex $(i_1,\ldots,i_k,j_1,\ldots,j_l)$ (resp.{} $(i_1,\ldots,i_k,i)$).

We conclude this section by briefly reviewing those aspects of the
Skinner-Rusk formalism for mechanics \cite{s83,sr83,sr83b} that survive in our
generalization to higher order field theory.

Let $Q$ be an $m$-dimensional smooth manifold and $q^{1},\ldots,q^{m}$ coordinates on it. Let $L\in C^{\infty}(TQ)$ be a Lagrangian function. Consider the
induced bundle $\tau_{0}^{\dag}:=\tau_{Q}^{\circ}(\tau_{Q}^{\ast}):T^{\dag
}:=\tau_{Q}^{\circ}(T^{\ast}Q)\longrightarrow TQ$ from the cotangent bundle
$\tau_{Q}^{\ast}:T^{\ast}Q\longrightarrow Q$ to $Q$, via the tangent bundle
$\tau_{Q}:TQ\longrightarrow Q$. Let $q:T^{\dag}\longrightarrow T^{\ast}Q$
be the canonical projection (see Diagram (\ref{DiagSkinRusk}))
\begin{equation}%
\begin{array}
[c]{c}%
\xymatrix{T^\dag \ar[d]_-{\tau_0^\dag} \ar[r]^-{q} & T^\ast Q \ar[d]^-{\tau^\ast_Q} \\
           TQ                \ar[r]^-{\tau_Q} & Q}
\end{array}
. \label{DiagSkinRusk}%
\end{equation}
On $T^{\dag}$ there is a canonical function $h\in C^{\infty}(T^{\dag})$
defined by $h(v,p):=p(v)$, $v\in T_{q}Q$, $p\in T_{q}^{\ast}Q$, $q\in Q$.
Consider also the function $E_{L}:=h-\tau_{0}^{\dag}{}^\ast(L)\in C^{\infty}(T^{\dag
})$. $E_{L}$ is locally given by $E_{L}:=p_{i}\dot{q}^{i}-L$, where
$\ldots,q^{i},\ldots,\dot{q}^{i},\ldots,p_{i},\ldots$ are standard coordinates
on $T^{\dag}$. Finally, put $\omega:=q^{\ast}(\omega_{0})\in\Lambda
^{2}(T^{\dag})$, $\omega_{0}\in\Lambda^{2}(T^{\ast}Q)$ being the canonical
symplectic form on $T^{\ast}Q$, which is locally given by $\omega_{0}%
=dp_{i}dq^{i}$. $\omega$ is a presymplectic form on $T^{\dag}$ whose kernel is
made of vector fields over $T^{\dag}$ which are vertical with respect to the
projection $q$. In the following, denote by $\mathbb{I}\subset\mathbb{R}$ a generic
open interval. For a curve $\gamma:\mathbb{I}\ni t\longmapsto\gamma(t)\in T^{\dag}$,
consider equations
\begin{equation}
i_{\dot{\gamma}}\gamma^\circ (\omega)+ \gamma^\circ(dE_{L})=0, \label{EqELH}%
\end{equation}
where $\dot{\gamma}\in\gamma^{\circ}(\mathrm{D}(T^{\dag}))$ is the tangent
field to $\gamma$. Equations (\ref{EqELH}) read locally
\[
\left\{
\begin{array}
[c]{l}%
\tfrac{d}{dt}q^{i}=\dot{q}^{i}\\
p_{i}=\tfrac{\partial L}{\partial\dot{q}^{i}}\\
\tfrac{d}{dt}p_{i}=\tfrac{\partial L}{\partial q^{i}}%
\end{array}
\right.  .
\]
In particular, for any solution $\gamma$ of Equations (\ref{EqELH}) as above,
$\tau_{Q}\circ\tau_{0}^{\dag}\circ\gamma:\mathbb{I}\longrightarrow Q$ is a solution of
the Euler-Lagrange equations determined by $L$. Notice that solutions of
Equations (\ref{EqELH}) can only take values in the submanifold
$\mathscr{P}\subset T^{\dag}$ defined as
\[
\mathscr{P}:=\{P\in T^{\dag}:\text{there exists }\Xi\in T_{P}T^{\dag}\text{
such that }i_{\Xi}\omega_{P}+(dE_{L})_{P}=0\},
\]
and that $\mathscr{P}$ is nothing but the graph of the Legendre transform
$FL:TQ\longrightarrow T^{\ast}Q$. Finally, consider $\mathscr{P}_{0}%
:=q(\mathscr{P})\subset T^{\ast}Q$. If $\mathscr{P}_{0}\subset T^{\ast}Q$
is a submanifold and $q:\mathscr{P}\longrightarrow\mathscr{P}_{0}$ a
submersion with connected fibers, then there exists a (unique) function $H\in
C^{\infty}(\mathscr{P}_{0})$ such that $q^{\ast}(H)=E_{L}|_{\mathscr{P}}$. Thus, for
a curve $\sigma:\mathbb{I}\ni t\longmapsto\sigma(t)\in\mathscr{P}_{0}$, we can consider
equations
\begin{equation}
i_{\dot{\sigma}}\sigma^\circ (\omega_{0})+\sigma^\circ(dH)=0, \label{EqH}%
\end{equation}
where $\dot{\sigma}\in\sigma^{\circ}(\mathrm{D}(T^{\dag}))$ is the tangent
field to $\sigma$. For any solution $\gamma:\mathbb{I}\longrightarrow Q$ of the
Euler-Lagrange equations, $FL\circ\dot{\gamma}:\mathbb{I}\longrightarrow\mathscr{P}_{0}%
$ is a solution of Equations (\ref{EqH}). If the map $q: \mathscr{P} \longrightarrow T^\ast Q$ has maximum rank (which happens iff the matrix $\left\Vert \partial^{2}L/\partial\dot
{q}^{i}\partial\dot{q}^{j}\right\Vert _{i}^{j}$ has maximum rank, i.e., $FL$
is a local diffeomorphism), then $\mathscr{P}_0 \subset T^\ast Q$ is an open submanifold, $H$ is a local function on $T^{\ast}Q$, and
Equations (\ref{EqH}) read locally
\[
\left\{
\begin{array}
[c]{l}%
\tfrac{d}{dt}q^{i}=\tfrac{\partial H}{\partial p_{i}}\\
\tfrac{d}{dt}p_{i}=-\tfrac{\partial H}{\partial q^{i}}%
\end{array}
\right.  ,
\]
which are Hamilton equations. In this case, for any solution $\sigma
:\mathbb{I}\longrightarrow T^{\ast}Q$ of Equations (\ref{EqH}), $\tau_{Q}^{\ast}%
\circ\sigma:\mathbb{I}\longrightarrow Q$ is a solution of the Euler-Lagrange equations.

\section{Geometry of Differential Equations\label{SecGeomPDE}}

In this section we recall basic facts about the geometric theory of PDEs. For more details see, for instance, \cite{b...99}.

Let $\pi:E\longrightarrow M$ be a fiber bundle, $\dim M=n$, $\dim E=m+n$. In the following we denote by $U \subset M$ a generic open subset. For
$0\leq l\leq k\leq\infty$, let $\pi_{k}:J^{k}\pi\longrightarrow M$ be
the bundle of $k$-jets of local sections of $\pi$, and $\pi_{k,l}:J^{k}%
\pi\longrightarrow J^{l}\pi$ the canonical projection. For any local section
$s:U\longrightarrow E$ of $\pi$, we denote
by $j_{k}s:U\longrightarrow J^{k}\pi$ its $k$th jet prolongation. For $x\in
U$, put $[s]_{x}^{k}:=(j_{k}s)(x)$. Any system of adapted to $\pi$ coordinates
$(\ldots,x^{i},\ldots,u^{\alpha},\ldots)$ on an open subset $V$ of $E$ gives
rise to a system of jet coordinates on $\pi_{k,0}^{-1}(V)\subset J^{k}\pi$
which we denote by $(\ldots,x^{i},\ldots,u^{\alpha}{}_{|I},\ldots)$ or simply
$(\ldots,x^{i},\ldots,u_{I}^{\alpha},\ldots)$ if this does not lead to
confusion, $|I|{}\leq k$, where we put $u_{\mathsf{O}%
}^{\alpha}:=u^{\alpha}$, $\alpha=1,\ldots,m$.

Now, let $k<\infty$, $\tau_{0}:T_{0}\longrightarrow J^{k}\pi$ be a vector
bundle, and $(\ldots,x^{i},\ldots,u_{I}^{\alpha},\ldots,v^{a},\ldots)$
adapted to $\tau_{0}$, local coordinates on $T_{0}$. A (possibly non-linear)
\emph{differential operator of order }$\leq k$\emph{ `acting on local sections
of }$\pi$\emph{, with values in }$\tau_{0}$\emph{'} (in short \emph{`from
}$\pi$\emph{ to }$\tau_{0}$\emph{'}) is a section $\Phi:J^{k}\pi
\longrightarrow T_{0}$ of $\tau_{0}$.

Let $\pi^{\prime}:E^{\prime}\longrightarrow M$ be another fiber bundle and
$\varphi:E\longrightarrow E^{\prime}$ a morphism of bundles. For any local
section $s:U\longrightarrow E$ of $\pi$,
$\varphi\circ s:U\longrightarrow E^{\prime}$ is a local section of
$\pi^{\prime}$. Therefore, for all $0\leq k\leq\infty$, $\varphi$ induces
a morphism $j_{k}\varphi:J^{k}\pi\longrightarrow J^{k}\pi^{\prime}$ of the
bundles $\pi_{k}$ and $\pi_{k}^{\prime}$ defined by $(j_{k}\varphi)[s]_{x}%
^{k}:=[\varphi\circ s]_{x}^{k}$, $x\in U$. Diagram
\[
\xymatrix{ J^l \pi \ar[r]^-{j_l \varphi} \ar[d]_-{\pi_{l,k}} & J^l \pi^\prime \ar[d]^-{\pi^\prime_{l,k}} \\
            J^k \pi \ar[r]^-{j_k \varphi}  & J^k \pi^\prime }
\]
commutes for all $0\leq k\leq l\leq\infty$. $j_{k}\varphi$ is called the
$k$\emph{th prolongation of }$\varphi$.

The above construction generalizes to differential operators as shown, for instance, in \cite{v09b}. If $\Phi$ is a differential operator as above, we denote by $\Phi^{(l)}$ its $l$th prolongation. Moreover, put $\mathscr{E}_{\Phi}:=\{\theta\in J^{k}\pi:\Phi
(\theta)=0\}$. $\mathscr{E}_{\Phi}$ is called the \emph{(system of) PDE(s)}
determined by $\Phi$. For $0\leq l\leq\infty$ put also $\mathscr{E}_{\Phi
}^{(l)}:=\mathscr{E}_{\Phi^{(l)}}\subset J^{k+l}\pi$. If $\mathscr{E}_{\Phi
}$ is locally defined by
\begin{equation}
\Phi^{a}(\ldots,x^{i},\ldots,u_{I}^{\alpha},\ldots)=0,\quad
a=1,\ldots,p
\end{equation}
$\ldots,\Phi^{a}:=\Phi^{\ast}(v^{a}),\ldots$ being local functions on
$J^{k}\pi$, then
$\mathscr{E}_{\Phi
}^{(l)}$ is locally defined by
\begin{equation}
(D_{J}\Phi^{a})(\ldots,x^{i},\ldots,u_{I}^{\alpha},\ldots)=0,\quad
a=1,\ldots,p,\;|J|{}\leq l, \label{Prolong}%
\end{equation}
where  $D_{(j_{1},\ldots j_{l})}:=D_{j_{1}}\circ\cdots\circ D_{j_{l}}$, and
$D_{j}:=\partial/\partial x^{j}+u_{Ij}^{\alpha}\partial/\partial u_{I}%
^{\alpha}$ is the $j$\emph{th total derivative, }$j,j_{1},\ldots
,j_{l}=1,\ldots,m$. $\mathscr{E}_{\Phi}^{(l)}$ is called the $l$th
prolongation of the PDE $\mathscr{E}_{\Phi}$. In the following we put
$\partial_{\alpha}^{I}:=\partial/\partial u_{I}^{\alpha}$, $\alpha=1,\ldots
,m$.

A local section $s$ of $\pi$ is a \emph{(local) solution of }%
$\mathscr{E}_{\Phi}$ iff, by definition, $\operatorname{im}j_{k}%
s\subset\mathscr{E}_{\Phi}$ or, which is the same, $\operatorname{im}%
j_{k+l}s\subset\mathscr{E}_{\Phi}^{(l)}$ for some $l\leq\infty$. Notice that
the $\infty$th prolongation of $\mathscr{E}_{\Phi}$, $\mathscr{E}_{\Phi
}^{(\infty)}\subset J^{\infty}\pi$, is an inverse limit of the sequence of
maps
\begin{equation}
\xymatrix@C=45pt{M  &\mathscr{E}_\Phi \ar[l]_-{\pi_k} & \cdots \ar[l]& \mathscr{E}^{(l)}_\Phi \ar[l]_-{\pi_{k+l,k+l-1}} & \mathscr{E}^{(l+1)}_\Phi \ar[l]_-{\pi_{k+l+1,k+l}} & \cdots \ar[l]
} \label{Einfty}%
\end{equation}
and consists of \textquotedblleft formal solutions\textquotedblright\ of
$\mathscr{E}_{\Phi}$, i.e., possibly non-converging Taylor series fulfilling
(\ref{Prolong}) for every $l$.

$J^{\infty}\pi$ is not a finite dimensional smooth manifold. However, it is a
\emph{pro-finite dimensional smooth manifold}. For an introduction to the
geometry of pro-finite dimensional smooth manifolds see \cite{v01} (see also
\cite{v09b}, and \cite{s89,t91} for different approaches). In the following we
will only consider \emph{regular }PDEs, i.e., PDEs $\mathscr{E}_{\Phi}$ such
that $\mathscr{E}_{\Phi}^{(\infty)}\subset J^{\infty}\pi$ is a smooth
pro-finite dimensional submanifold in $J^{\infty}\pi$, i.e., $\pi_{\infty
,l}(\mathscr{E}_{\Phi}^{(\infty)})\subset J^{l}\pi$ is a smooth submanifold
and $\pi_{l+1,l}:\pi_{\infty,l+1}(\mathscr{E}_{\Phi}^{(\infty)}%
)\longrightarrow\pi_{\infty,l}(\mathscr{E}_{\Phi}^{(\infty)})$ is a smooth bundle for all $l \geq 0$.

There is a dual concept to the one of a pro-finite dimensional manifold, i.e.,
the concept of a \emph{filtered smooth manifold} which will be used in the
following. We do not give here a complete definition of a filtered manifold,
which would take too much space. Rather, we will just outline it. Basically, a
filtered smooth manifold is a(n equivalence class of) set(s) $\mathscr{O}$
together with a sequence of embeddings of closed submanifolds%
\begin{equation}
\xymatrix@C=40pt{\ \mathscr{O}_0\  \ar@{^{(}->}[r]^-{i_{0,1}} &\ \mathscr{O}_1\  \ar@{^{(}->}[r]^-{i_{0,1}} &\ \cdots\ \ar@{^{(}->}[r]& \ \mathscr{O}_{k-1}\  \ar@{^{(}->}[r]^-{i_{k-1,k}} &\ \mathscr{O}_{k}\  \ar@{^{(}->}[r]^-{i_{k,k+1}} &\ \cdots 
} \label{Ofilt}%
\end{equation}
and inclusions $i_{k}:\mathscr{O}_{k}\hookrightarrow\mathscr{O}$, $k\geq0$,
such that $\mathscr{O}$ (together with the $i_{k}$'s) is a direct limit of
(\ref{Ofilt}). The tower of
algebra epimorphisms
\begin{equation}
\xymatrix@C=40pt{C^\infty(\mathscr{O}_{0}) & \cdots \ar[l] & C^\infty(\mathscr{O}_{k}) \ar[l]_-{i_{k-1,k}^\ast} & C^\infty(\mathscr{O}_{k+1}) \ar[l]_-{i_{k,k+1}^\ast} & \cdots \ar[l]
} \label{COfilt}%
\end{equation}
is associated to sequence (\ref{Ofilt}).
We define $C^{\infty}(\mathscr{O})$ to be the inverse limit of the tower
(\ref{COfilt}). Every element in $C^{\infty}(\mathscr{O})$ is naturally a
function on $\mathscr{O}$. Thus, we interpret $C^{\infty}(\mathscr{O})$ as the
algebra of \emph{smooth functions on }$\mathscr{O}$. Clearly, there are
canonical \textquotedblleft restriction homomorphisms\textquotedblright%
\ $i_{k}^{\ast}:C^{\infty}(\mathscr{O})\longrightarrow C^{\infty
}(\mathscr{O}_{k})$, $k\geq0$. Differential calculus over $\mathscr{O}$ may
then be introduced as \emph{differential calculus over }$C^{\infty
}(\mathscr{O})$ \cite{v01} respecting the sequence (\ref{COfilt}). Since the
main constructions (smooth maps, vector fields, differential forms, jets and
differential operators, etc.) of such calculus and their properties do not
look very different from the analogous ones in finite-dimensional differential
geometry we will not insist on this. Just as an instance, we report here the
definition of a differential form $\omega$ on $\mathscr{O}$: it is just a
sequence of differential forms $\omega_{k}\in\Lambda(\mathscr{O}_{k})$,
$k\geq0$, such that $i_{k-1,k}^{\ast}(\omega_{k})=\omega_{k-1}$ for all $k$.

Finally, notice that, allowing for the $\mathscr{O}_{k}$'s in (\ref{Ofilt}) to
be pro-finite dimensional manifolds, we obtain a more general object than both
a pro-finite dimensional and a filtered manifold. We will generically refer to such an object
as an \emph{infinite dimensional smooth manifold} or even just
a \emph{smooth manifold} if this does not lead to confusion. Our main example of
such a kind of infinite dimensional manifold will be presented in the
beginning of Section \ref{SecLH}.

\section{The Cartan Distribution and the Lagrangian
Formalism\label{SecCartLagr}}

Let $\pi:E\longrightarrow M$ and $\Phi$ be as in the previous section. In the
following we will simply write $J^{k}$ for $J^{k}\pi$, $k \leq \infty$, and
$\mathscr{E}$ for $\mathscr{E}_{\Phi}^{(\infty)}$. $\mathscr{E}$ will be
referred to simply as a PDE (imposed on sections of $\pi$) if this does not
lead to confusion. Notice that for $\Phi=0$, $\mathscr{E}=\mathscr{E}_{\Phi
}^{(\infty)}=J^{\infty}$.

Recall that $J^{\infty}$ is canonically endowed with the Cartan distribution \cite{b...99}
\[
\mathscr{C}:J^{\infty}\ni\theta\longmapsto\mathscr{C}_{\theta}\subset
T_{\theta}J^{\infty}%
\]
which is locally spanned by total derivatives, $D_{i}$, $i=1,\ldots,n$.
$\mathscr{C}$ is a flat connection in $\pi_{\infty}$ which we call the
\emph{Cartan connection}. Moreover, it restricts to $\mathscr{E}$ in the sense
that $\mathscr{C}_{\theta}\subset T_{\theta}\mathscr{E}$ for any $\theta
\in\mathscr{E}$. Therefore, the (infinite prolongation of) any PDE is
naturally endowed with an involutive distribution whose $n$-dimensional
integral submanifolds are of the form $j_{\infty}s$, with $s:U\longrightarrow
E$ a (local) solution of $\mathscr{E}_{\Phi}$. In
the following we will identify the space of $n$-dimensional integral
submanifolds of $\mathscr{C}$ and the space of local solutions of
$\mathscr{E}_{\Phi}$.

Let $\pi^{\prime}:E^{\prime}\longrightarrow M$ be another bundle and
$\mathscr{E}^{\prime}\subset J^{\infty}\pi^{\prime}$ (the infinite
prolongation of) a PDE imposed on sections of $\pi^{\prime}$. A smooth map
$F:\mathscr{E}^{\prime}\longrightarrow\mathscr{E}$ is called a \emph{morphism
of PDEs} iff it respects the Cartan distributions, i.e., $(d_{\theta^{\prime}%
}F)(\mathscr{C}_{\theta^{\prime}}) = \mathscr{C}_{F(\theta^{\prime})}$
for any $\theta^{\prime}\in\mathscr{E}^{\prime}$. The idea of \emph{non-local
variables} in the theory of PDEs can be formalized geometrically by special
morphisms of PDEs called \emph{coverings} \cite{kv89} (see also \cite{i06}). A
\emph{covering} is a morphism $\psi:\widehat{\mathscr{E}}\longrightarrow
\mathscr{E}$ of PDEs which is surjective and submersive. A covering
$\psi:\widehat{\mathscr{E}}\longrightarrow\mathscr{E}$ clearly sends local
solutions of $\widehat{\mathscr{E}}$ to local solutions of $\mathscr{E}$. If
there exists a covering $\psi:\widehat{\mathscr{E}}\longrightarrow\mathscr{E}$
of PDEs we also say that the PDE $\widehat{\mathscr{E}}$\emph{ covers the PDE
}$\mathscr{E}$\emph{ (via }$\psi$\emph{)}. Fiber coordinates on the total
space $\widehat{\mathscr{E}}$ of a covering $\psi:\widehat{\mathscr{E}}%
\longrightarrow\mathscr{E}$ are naturally interpreted as non-local variables
on $\mathscr{E}$. Also notice that given a solution $s$ of the PDE
$\mathscr{E}$, a covering $\psi:\widehat{\mathscr{E}}\longrightarrow
\mathscr{E}$ determines a whole family of solutions of $\widehat{\mathscr{E}}$
\textquotedblleft projecting onto $s$ via $\psi$\textquotedblright, so that
$\psi$ may be interpreted, to some extent, as a fibration over the space of
solutions of $\mathscr{E}$. Many relevant constructions in the theory of PDEs (including Lax pairs,
B\"{a}cklund transformations, etc.) are duly formalized in geometrical terms
by using coverings.

The Cartan distribution and the fibered structure $\pi_{\infty}:J^{\infty
}\longrightarrow M$ of $J^{\infty}$ determine a splitting of the tangent
bundle $TJ^{\infty}\longrightarrow J^{\infty}$ into the Cartan or horizontal
part $\mathscr{C}$ and the vertical (with respect to $\pi_{\infty}$) part.
Accordingly, the de Rham complex of $J^{\infty}$,
$(\Lambda(J^{\infty}),d)$, splits in the \emph{variational bi-complex} $(\mathscr{C}^{\bullet
}\Lambda\otimes\overline{\Lambda}{},\overline{d},d^{V})$, (here and in
what follows tensor products will be always over $C^{\infty}(J^{\infty})$ if
not otherwise specified), where $\mathscr{C}^{\bullet}\Lambda$ and $\overline{\Lambda}{}^\bullet$ are the algebras of \emph{Cartan forms} and \emph{horizontal forms} respectively. $d^{V}$ and $\overline{d}$ are the \emph{vertical}
and the \emph{horizontal de Rham differential}, respectively (see, for instance, \cite{b...99} for details). The variational bicomplex allows a cohomological formulation of the calculus of
variations \cite{v84,b...99,v01,a92}. In the second part of this section we
briefly review it. 

In the following we will understand isomorphism $\Lambda(J^{\infty}%
)\simeq\mathscr{C}^{\bullet}\Lambda\otimes\overline{\Lambda}{}$. The complex
\[
\xymatrix{0 \ar[r] &  C^\infty(J^\infty) \ar[r]^-{\overline{d}} & \overline{\Lambda}{}^1 \ar[r]^-{\overline{d}} & \cdots \ar[r] & \overline{\Lambda}{}^q \ar[r]^-{\overline{d}} & \overline{\Lambda}{}^{q+1} \ar[r]^-{\overline{d}} & \cdots}
\]
is called the \emph{horizontal de Rham complex}. An element $\mathscr{L}\in
\overline{\Lambda}{}^{n}$ is naturally interpreted as a \emph{Lagrangian
density} and its cohomology class $[\mathscr{L}]\in \overline{H}{}^{n}:=H^{n}(\overline{\Lambda},\overline{d})$ as an \emph{action functional}
on sections of $\pi$. The associated Euler-Lagrange equations can then be
obtained as follows.

Consider the complex
\begin{equation}
\xymatrix{0 \ar[r] &  \mathscr{C}\Lambda^1 \ar[r]^-{\overline{d}} & \mathscr{C}\Lambda^1\otimes\overline{\Lambda}{}^1 \ar[r]^-{\overline{d}} & \cdots \ar[r] & \mathscr{C}\Lambda^1\otimes\overline{\Lambda}{}^q \ar[r]^-{\overline{d}} & \cdots}, \label{ComplCL1}%
\end{equation}
and the $C^{\infty}(J^{\infty})$-submodule $\varkappa^{\dag}\subset
\mathscr{C}\Lambda^{1}\otimes\overline{\Lambda}{}^{n}$ generated by elements
in $\mathscr{C}\Lambda^{1}\otimes\overline{\Lambda}{}^{n}\cap\Lambda(J^{1}\pi)$. $\varkappa^{\dag}$ is locally spanned by elements
$(du^{\alpha}-u_{i}^{\alpha}dx^{i})\otimes d^{n}x$, where we put
$d^{n}x:=dx^{1}\cdots dx^{n}$.

\begin{thm}
\label{Theorem0}\cite{v84} Complex (\ref{ComplCL1}) is acyclic in the $q$th
term, for $q\neq n$. Moreover, for any $\omega\in\mathscr{C}\Lambda^{1}%
\otimes\overline{\Lambda}{}^{n}$ there exists a unique element $\boldsymbol{E}%
_{\omega}\in\varkappa^{\dag}\subset\mathscr{C}\Lambda^{1}\otimes
\overline{\Lambda}{}^{n}$ such that $\boldsymbol{E}_{\omega}-\omega
=\overline{d}\vartheta$ for some $\vartheta\in\mathscr{C}\Lambda^{1}%
\otimes\overline{\Lambda}{}^{n-1}$ and the correspondence $H^{n}%
(\mathscr{C}\Lambda^{1}\otimes\overline{\Lambda},\overline{d})\ni\lbrack
\omega]\longmapsto\boldsymbol{E}_{\omega}\in\varkappa^{\dag}$ is a vector
space isomorphism. In particular, for $\omega=d^{V}\!\mathscr{L}$,
$\mathscr{L}\in\overline{\Lambda}{}^{n}$ being a Lagrangian density locally
given by $\mathscr{L}=Ld^{n}x$, $L$ a local function on $C^{\infty}(J^{\infty
})$, $\boldsymbol{E}(\mathscr{L}):=\boldsymbol{E}_{\omega}$ is locally given
by $\boldsymbol{E}(\mathscr{L})=\tfrac{\delta L}{\delta u^{\alpha}}%
(du^{\alpha}-u_{i}^{\alpha}dx^{i})\otimes d^{n}x$ where $\tfrac{\delta
L}{\delta u^{\alpha}}:=(-)^{|I|}D_{I}\partial_{\alpha}^{I}L$ are the
Euler-Lagrange derivatives of $L$.
\end{thm}

In view of the above theorem, $\boldsymbol{E}(\mathscr{L})$ does not depend on
the choice of $\mathscr{L}$ in a cohomology class $[\mathscr{L}]\in$
$\overline{H}{}^{n}$ and it is naturally interpreted as the left hand side of
the Euler-Lagrange (EL) equations determined by $\mathscr{L}$. In the
following we will denote by $\mathscr{E}_{EL}\subset J^{\infty}$ the (infinite
prolongation of the) EL equations determined by a Lagrangian density. Any
$\vartheta\in\mathscr{C}\Lambda^{1}\otimes\overline{\Lambda}{}^{n-1}$ such
that
\begin{equation}
\boldsymbol{E}(\mathscr{L})-d^{V}\!\mathscr{L}=\overline{d}\vartheta
\label{Legendre}%
\end{equation}
will be called a \emph{Legendre form} \cite{av04}. Equation (\ref{Legendre})
may be interpreted as the \emph{first variation formula} for the Lagrangian
density $\mathscr{L}$. In this respect, the existence of a global Legendre
form was first discussed in \cite{k80}.

\begin{rem}
\label{Remark2}Notice that, if $\vartheta\in\mathscr{C}\Lambda^{1}%
\otimes\overline{\Lambda}{}^{n-1}$ is a Legendre form for a Lagrangian density
$\mathscr{L}\in\overline{\Lambda}{}^{n}$, then $\vartheta+d^{V}\!\varrho$ is a
Legendre form for the $\overline{d}$-cohomologous Lagrangian density $\mathscr{L}+\overline
{d}\varrho$, $\varrho\in\overline{\Lambda}{}^{n-1}$, which determines the same
EL equations as $\mathscr{L}$. Moreover, any two Legendre forms $\vartheta
,\vartheta^{\prime}$ for the same Lagrangian density differ by a $\overline
{d}$-closed, and, therefore, $\overline{d}$-exact form, i.e., $\vartheta
-\vartheta^{\prime}=\overline{d}\lambda$, for some $\lambda\in
\mathscr{C}\Lambda^{1}\otimes\overline{\Lambda}{}^{n-2}$.
\end{rem}

\begin{rem}
\label{Remark1}Finally, notice that complex (\ref{ComplCL1}) restricts to
\emph{holonomic sections} $j_{\infty}s$ of $\pi_{\infty}$, $s$ being a local
sections of $\pi$, in the sense that, for any such $s$, there is a (unique)
complex
\begin{equation}
\xymatrix{0 \ar[r] &  \mathscr{C}\Lambda^1|_j \ar[r]^-{\overline{d}|_j} & \mathscr{C}\Lambda^1\otimes\overline{\Lambda}{}^1|_j \ar[r]^-{\overline{d}|_j} & \cdots \ar[r] & \mathscr{C}\Lambda^1\otimes\overline{\Lambda}{}^q|_j \ar[r]^-{\overline{d}|_j} & \cdots}, \label{ComplCL1j}%
\end{equation}
where $j:=j_{\infty}s$, such that the restriction map $\mathscr{C}\Lambda
^{1}\otimes\overline{\Lambda}{}\longrightarrow\mathscr{C}\Lambda^{1}%
\otimes\overline{\Lambda}{}|_{j}\simeq \mathscr{C}\Lambda^{1}|_{j}\otimes
_{C^{\infty}(M)}\Lambda(M)$ is a morphism of complexes. Moreover, complex
(\ref{ComplCL1j}) is acyclic in the $q$th term and the correspondence defined
by $H^{n}(\mathscr{C}\Lambda^{1}\otimes\overline{\Lambda}{}^{n}|_{j}%
,\overline{d}|_{j})\ni\lbrack\omega|_{j}]\longmapsto\boldsymbol{E}_{\omega
}|_{j}\in\varkappa^{\dag}|_{j}$, $\omega\in\mathscr{C}\Lambda^{1}%
\otimes\overline{\Lambda}{}^{n}$, is a vector space isomorphism.
\end{rem}

\section{\label{SecPDHam}Partial Differential Hamiltonian Systems}

In \cite{v09} we defined a PD analogue
of the concept of Hamiltonian system on an abstract symplectic manifold which
we called a \emph{PD-Hamiltonian system}. In this section we briefly review
those definitions and results in \cite{v09} which we will need in the following.

Let $\alpha:P\longrightarrow M$ be a fiber bundle, $A:=C^{\infty}(P)$,
$x^{1},\ldots,x^{n}$ coordinates on $M$, $\dim M=n$, and $q^{1},\ldots,q^{m}$
fiber coordinates on $P$, $\dim P=n+m$. Denote by $C(P,\alpha)$ the space of
(Ehresmann) connections in $\alpha$. $C(P,\alpha)$ identifies canonically with
the space of sections of the first jet bundle $\alpha_{1,0}:J^{1}%
\alpha\longrightarrow P$ and in the following we will understand such
identification. In particular, for $\nabla\in C(P,\alpha)$, we put
$\ldots,\nabla_{i}^{A}:=\nabla^{\ast}(q_{i}^{A}),\ldots$, $\ldots,q_{i}%
^{A},\ldots$ being jet coordinates on $J^{1}\alpha$.

Denote by $\Lambda_{1}=\bigoplus_{k}\Lambda_{1}^{k}\subset\Lambda(P)$ the
differential (graded) ideal in $\Lambda(P)$ made of differential forms on $P$
vanishing when pulled-back to fibers of $\alpha$, by $\Lambda_{p}%
=\bigoplus_{k}\Lambda_{p}^{k}$ its $p$-th exterior power, $p\geq0$, and by
$V\!\Lambda(P,\alpha)=\bigoplus_{k}V\!\Lambda^{k}(P,\alpha)$ the quotient
differential algebra $\Lambda(P)/\Lambda_{1}$, $d^{V}:V\!\Lambda
(P,\alpha)\longrightarrow V\!\Lambda(P,\alpha)$ being its (quotient) differential.

\begin{rem}
\label{Remark3}For instance, if $\alpha=\pi_{\infty}:P=J^{\infty
}\longrightarrow M$, then, using the Cartan connection $\mathscr{C}\in
C(J^{\infty},\pi_{\infty})$, one can canonically identify $V\!\Lambda
^{1}(J^{\infty},\pi_{\infty})$ with $\mathscr{C}\Lambda^{1}$ and $d^{V}$ with
the vertical de Rham differential. More generally, for any $k\geq0$,
$V\!\Lambda^{1}(J^{k},\pi_{k})\otimes_{C^{\infty}(J^{k}\pi)}C^{\infty}%
(J^{k+1}\pi)$ identifies canonically with the $C^{\infty}(J^{k+1}\pi)$-module
$\mathscr{C}\Lambda^{1}\cap\Lambda(J^{k+1}\pi)$ of $(k+1)$th order
\emph{Cartan forms}.
\end{rem}

Now, for any $k\geq0$, put $\Omega^k(P,\alpha) := \Lambda^{k+n-1}_{n-1}$ and $\underline{\Omega}^k(P,\alpha) := \Omega^k(P,\alpha)/ \Lambda_n^{k+n-1}$. It is easy to show that $\underline{\Omega}^k(P,\alpha) \simeq V\!\Lambda^k(P,\alpha)\otimes_A \Lambda^{n-1}_{n-1}$. An element $\omega \in \Omega^k (P,\alpha)$ determines an affine map
\begin{equation}
C(P,\alpha) \ni \nabla \longmapsto i_\nabla \omega := \mathfrak{p}_{\nabla}(\omega)\in V\!\Lambda
^{k-1}(P,\alpha)\otimes_{A}\Lambda_{n}^{n}, \label{AffMap}
\end{equation} 
where
\[\mathfrak{p}_{\nabla}:\Lambda(P)\longrightarrow V\!\Lambda^{k-1}(P,\alpha)\otimes
_{A}\Lambda_{n}^{n}\] is the canonical projection determined by the
connection $\nabla$. The linear part of the affine map (\ref{AffMap}) naturally identify with the class $\omega + \Lambda_n^{k+n-1}$ in $\underline{\Omega}^k(P,\alpha)$ (see \cite{v09} for details).
Notice that, since (\ref{AffMap}) is affine, it is actually point-wise and, therefore, can be restricted to maps.
Namely, if $F:P_{1}\longrightarrow P$ is a smooth map, $\square\in F^{\circ
}(C(P,\alpha))$, then an element $i_{\square}F^{\circ
}(\omega)\in F^{\circ}(V\!\Lambda^{k-1}(P,\alpha)\otimes_{A}\Lambda
_{n}^{n})$ is defined in an obvious way.

\begin{defn}
A \emph{PD-Hamiltonian system} on the fiber bundle $\alpha:P\longrightarrow M$
is an element $\omega\in\Omega^{2}(P,\alpha)$ such that $d\omega = 0$. The first order
PDEs
\[
i_{j_{1}\sigma}\omega|_{\sigma}=0
\]
on (local) sections $\sigma$ of $\alpha$ are called the \emph{PD-Hamilton
equations determined by }$\omega$. Geometrically, they correspond to the
submanifold
\[
\mathscr{E}_{\omega}^{(0)}:=\{\theta\in J^{1}\alpha:i_{\theta}\omega
_{p}=0\text{, }p=\alpha_{1,0}(\theta)\}\subset J^{1}\alpha.
\]

\end{defn}

Let $\omega$ be a PD-Hamiltonian system on the bundle
$\alpha:P\longrightarrow M$ and consider the subset $P_{1}:=\alpha
_{1,0}(\mathscr{E}_{\omega}^{(0)})\subset P$. In the following we will assume
$P_{1}\subset P$ to be a submanifold and $\alpha_{1}:=\alpha|_{P_{1}}%
:P_{1}\longrightarrow M$ to be a subbundle of $\alpha$. $\alpha_{1}$ is called
\emph{the first constraint subbundle of }$\omega$.

As an example, consider the following canonical constructions. Let
$\alpha:P\longrightarrow M$ be a fiber bundle and $\ldots,q^{A},\ldots$
fiber coordinates on $P$. $\Omega^{1}(P,\alpha)$ (resp. $\underline{\Omega}%
^{1}(P,\alpha)$) is the $C^{\infty}(P)$-module of sections of a vector bundle
$\mu_{0}\alpha:\mathscr{M}\alpha\longrightarrow P$ (resp. $\tau_{0}^{\dag}%
\alpha:J^{\dag}\alpha\longrightarrow P$), called the \emph{multimomentum bundle of}
$\alpha$ (resp. the \emph{reduced multimomentum bundle of }$\alpha$). Recall that
there is a tautological element $\Theta_{\alpha} \in \Omega^{1}(\mathscr{M}\alpha
,\mu\alpha)$ (resp. $\underline{\Theta}{}_{\alpha}\in\underline{\Omega}^{1}(J^{\dag
}\alpha,\tau^{\dag}\alpha)$), where $\mu\alpha:=\alpha\circ\mu_{0}\alpha$ (resp. $\tau^{\dag
}\alpha:=\alpha\circ\tau_{0}^{\dag}\alpha$), which
in standard coordinates $\ldots,x^{i},\ldots,q^{A},\ldots,p_{A}%
^{i},\ldots,p$ on $\mathscr{M}\alpha$ (resp. $\ldots,x^{i},\ldots,q^{A
},\ldots,p_{A}^{i},\ldots$ on $J^{\dag}\alpha$) is given by
\[
\Theta_{\alpha}=p_{A}^{i}dq^A d^{n-1}x_{i}-pd^{n}x\quad\text{(resp. }
\underline{\Theta}{}_{\alpha}=p_{A}^{i}d^{V}\!q^{A}\otimes d^{n-1}%
x_{i}\text{),}%
\]
where $d^{n-1}x_{i}:=i_{\partial/\partial x^{i}}d^{n}x$ \cite{gim98}. $d\Theta_{\alpha
}$ is a PD-Hamiltonian system on
$\mu\alpha$ locally given by
\[
d\Theta_{\alpha}=dp_{A}^{i}dq^{A}d^{n-1}x_{i}-dpd^{n}x.
\]
Notice that $d\Theta_{\alpha}$ determines empty PD-Hamilton equations.

\begin{exmp}
\label{Example1}A PD-Hamiltonian system is canonically determined, on the fiber
bundle $\alpha:P\longrightarrow M$, by the following data: a
connection $\nabla$ in $\alpha$ and a differential form
$\mathscr{L}\in\Lambda_{n}^{n}$. Let $\ldots,q^{A},\ldots$ be fiber
coordinates in $P$ and $\ldots,x^{i},\ldots,q^{A},\ldots,p_{A}^{i},\ldots,p$
(resp. $\ldots,x^{i},\ldots,q^{A},\ldots,p_{A}^{i},\ldots$) standard
coordinates in $\mathscr{M}\alpha\ $(resp. $J^{\dag}\alpha$). Let
$\mathscr{L}$ be locally given by $\mathscr{L}=Ld^{n}x$, $L$ a local
function on $P$. Obviously, $\nabla$ determines a section $\Sigma_\nabla : J^{\dag
}\alpha \longrightarrow \mathscr{M}\alpha$ of the projection $\mathscr{M}\alpha\longrightarrow J^{\dag
}\alpha$,
which in local standard coordinates reads $\Sigma_{\nabla}^{\ast}(p)=p_{A}%
^{i}\nabla_{i}^{A}$. Put $\Theta_{\nabla}:=\Sigma_{\nabla}^{\ast}%
(\Theta_{\alpha})$. In local
standard coordinates, $\Theta_{\nabla}=p_{A}^{i}dq^{A}d^{n-1}x_{i}-p_{A}%
^{i}\nabla_{i}^{A}d^{n}x.$ Put also,
\[
\Theta_{\mathscr{L},\nabla}:=\Theta_{\nabla}+(\tau_{0}^{\dag}\alpha)^{\ast
}(\mathscr{L}).
\]
Locally, $\Theta_{\mathscr{L},\nabla}=p_{A}^{i}dq^{A}d^{n-1}x_{i}%
-E_{\mathscr{L},\nabla}d^{n}x$, where $E_{\mathscr{L},\nabla}:=p_{A}^{i}%
\nabla_{i}^{A}-L$. Finally, consider $\omega_{\mathscr{L},\nabla}%
:=d\Theta_{\mathscr{L},\nabla}$. Locally,
\[
\omega_{\mathscr{L},\nabla}=dp_{A}^{i}dq^{A}d^{n-1}x_{i}%
-dE_{\mathscr{L},\nabla}d^{n}x.
\]
$\omega_{\mathscr{L},\nabla}$ is the \emph{PD-Hamiltonian system on }%
$\tau^{\dag}\alpha$ \emph{determined by }$\nabla$\emph{ and }$\mathscr{L}$.
The associated PD-Hamilton equations read locally
\[
\left\{
\begin{array}
[c]{l}%
p_{A}^{i},_{i}=\tfrac{\partial }{\partial q^A}L-p_{B}^{i}\tfrac{\partial }{\partial q^A}\nabla_{i}^{B}\\
q^{A},_{i}=\nabla_{i}^{A}%
\end{array}
\right.  ,
\]
where we denoted by \textquotedblleft${}\bullet{},_{i}$\textquotedblright\ 
the partial derivative of \textquotedblleft${}\bullet{}$\textquotedblright%
\ with respect to the $i$th independent variable $x^{i}$, $i=1,\ldots,n$.
\end{exmp}

We conclude this section by discussing two examples of morphisms of PDEs
coming from the theory of PD-Hamiltonian systems.

\begin{exmp}
Let $\alpha:P\longrightarrow M$ be a fiber bundle, $\omega$ a PD-Hamiltonian system on it, $\alpha^{\prime
}:P^{\prime}\longrightarrow M$ another fiber bundle, $\beta:P^{\prime
}\longrightarrow P$ a surjective, submersive, fiber bundle morphism, and
$\omega^{\prime}:=\beta^{\ast}(\omega)$. $\omega^{\prime}$ is a PD-Hamiltonian system on $\alpha^{\prime}$.
Denote by $\mathscr{E}\subset J^{\infty}\alpha$ (resp. $\mathscr{E}^{\prime
}\subset J^{\infty}\alpha^{\prime}$) the $\infty$th prolongation of the
PD-Hamilton equations determined by $\omega$ (resp. $\omega^{\prime}$). We
want to compare $\mathscr{E}$ and $\mathscr{E}^{\prime}$. In order to do this,
notice, preliminarily, that $J^{\infty}\alpha^{\prime}$ covers $J^{\infty
}\alpha$ via $j_{\infty}\beta:J^{\infty}\alpha^{\prime}\longrightarrow
J^{\infty}\alpha$. Moreover, it can be easily checked that a local section
$\sigma^{\prime}$ of $\alpha^{\prime}$ is a solution of $\mathscr{E}^{\prime}$
iff the section $\beta\circ\sigma^{\prime}$ of $\alpha$ is a solution of
$\mathscr{E}$. We now prove the formal version of this fact.

\begin{prop}
\label{Prop2}$(j_{\infty}\beta)(\mathscr{E}^{\prime})\subset\mathscr{E}$ and
$j_{\infty}\beta:\mathscr{E}^{\prime}\longrightarrow\mathscr{E}$ is a covering.
\end{prop}

\begin{pf}
Consider $j_{1}\beta:J^{1}\alpha^{\prime}\longrightarrow J^{1}\alpha$. It is
easy to check that $\mathscr{E}_{\omega^{\prime}}^{(0)}=(j_{1}\beta
)^{-1}(\mathscr{E}_{\omega}^{(0)})\subset J^{1}\alpha^{\prime}$. Similarly,
$\mathscr{E}^{\prime}=(j_{\infty}\beta)^{-1}(\mathscr{E})\subset J^{\infty
}\alpha^{\prime}$. In particular, $j_{\infty}\beta:\mathscr{E}^{\prime
}\longrightarrow\mathscr{E}$ is the \textquotedblleft
restriction\textquotedblright\ of $j_{\infty}\beta:J^{\infty}\alpha^{\prime
}\longrightarrow J^{\infty}\alpha$ to $\mathscr{E}\subset J^{\infty}\alpha$
and, therefore, is a covering.
\end{pf}
\end{exmp}

\begin{exmp}
Let $\alpha:P\longrightarrow M$, $\omega$ and
$\mathscr{E}\subset J^{\infty}\alpha$ be as in the above example, and
$\alpha_{1}:P_{1}\longrightarrow M$ the first constraint subbundle of $\omega
$. Put $\omega_{1}:=i_{P_{1}}^{\ast}(\omega)$. $\omega_{1}$ is a PD-Hamiltonian system on
$\alpha_{1}$. Denote by $\mathscr{E}_{1}\subset J^{\infty}\alpha_{1}$ the
$\infty$th prolongation of the PD-Hamilton equations determined by $\omega
_{1}$. We want to compare $\mathscr{E}$ and $\mathscr{E}_{1}$. In order to do
this, notice, preliminarily, that $J^{\infty}\alpha_{1}$ may be understood as
a submanifold in $J^{\infty}\alpha$ via $j_{\infty}i_{P_{1}}:J^{\infty}%
\alpha_{1}\hookrightarrow J^{\infty}\alpha$. Moreover, it can be easily
checked that any solution of $\mathscr{E}$ is also a solution of
$\mathscr{E}_{1}$ (while the vice-versa is generically untrue). We now prove
the formal version of this fact.

\begin{prop}
\label{Prop1}$\mathscr{E}\subset\mathscr{E}_{1}$.
\end{prop}

\begin{pf}
Recall that the projection $\alpha_{1,0}:J^{1}\alpha\longrightarrow P$ sends
$\mathscr{E}_{\omega}^{(0)}$ to $P_{1}$. As a consequence, $\mathscr{E}\subset
J^{\infty}\alpha_{1}$. Moreover, by definition of $\infty$th prolongation of a
PDE, it is easy to check that
\begin{align*}
\mathscr{E}  &  =\mathscr{E}\cap J^{\infty}\alpha_{1}\\
&  =\{\theta=[\sigma]_{x}^{\infty}\in J^{\infty}\alpha_{1}:\text{ }%
[i_{j_{1}\sigma}\omega|_{\sigma}]_{x}^{\infty}=0,\text{ }x\in M\}\\
&  \subset\{\theta=[\sigma]_{x}^{\infty}\in J^{\infty}\alpha_{1}:\text{
}[i_{j_{1}\sigma}\omega_{1}|_{\sigma}]_{x}^{\infty}=0,\text{ }x\in M\}\\
&  =\mathscr{E}_{1}.
\end{align*}

\end{pf}
\end{exmp}

\section{Lagrangian-Hamiltonian Formalism\label{SecLH}}

In this section we show that the Skinner-Rusk mixed Lagrangian-Hamiltonian
formalism for first order mechanics \cite{s83,sr83,sr83b} (see Section
\ref{SecConvSkinRusk}) is straightforwardly generalized to higher order
Lagrangian field theories.

First of all, let us present our main example of a filtered manifold. Let
$\pi:E\longrightarrow M$ be a fiber bundle. Consider the infinite jet
bundle $\pi_{\infty}:J^{\infty}\longrightarrow M$ for which $\Lambda_{q}%
^{q}=\overline{\Lambda}{}^{q}$, $q\geq0$. Moreover, the $C^{\infty}(J^{\infty
})$-module $\underline{\Omega}^{1}(J^{\infty},\pi_{\infty})\simeq \mathscr{C}\Lambda
^{1}\otimes\overline{\Lambda}{}^{n-1}$ is canonically filtered by vector
subspaces $W_{k}:=\mathscr{C}\Lambda^{1}\otimes\overline{\Lambda}{}^{n-1}%
\cap\Lambda(J^{k+1}\pi)$, $k\geq0$. Denote by $\underline{\Omega}_{k}%
^{1}\subset\underline{\Omega}^{1}(J^{\infty},\pi_{\infty})$ the $C^{\infty
}(J^{\infty})$-submodule generated by $W_{k}$, $k\geq0$. Then, for all $k$,
$\underline{\Omega}_{k}^{1}$ is canonically isomorphic to $C^{\infty
}(J^{\infty})\otimes_{C^{\infty}(J^{k+1})}W_{k}$ and
\begin{equation}
\underline{\Omega}_{0}^{1}\subset\underline{\Omega}_{1}^{1}\subset
\cdots\subset\underline{\Omega}_{k}^{1}\subset\underline{\Omega}_{k+1}%
^{1}\subset\cdots\subset\underline{\Omega}^{1}(J^{\infty},\pi_{\infty}),
\label{Eq5}%
\end{equation}
is a sequence of $C^{\infty}(J^{\infty})$-submodules. Notice that, for any
$k$, $\underline{\Omega}_{k}^{1}$ is the module of sections of a
finite-dimensional vector bundle $J_{k}^{\dag}%
\longrightarrow J^{\infty}$. Moreover, the inclusions (\ref{Eq5}) determine
inclusions
\[
J_{0}^{\dag}\subset J_{1}^{\dag}\subset\cdots\subset J_{k}^{\dag}\subset
J_{k+1}^{\dag}\subset\cdots
\]
of vector bundles. $J^{\dag}:=\bigcup_{k}J_{k}^{\dag}$ is then an infinite
dimensional (filtered) manifold and the canonical projection $\tau_{0}^{\dag
}:J^{\dag}\longrightarrow J^{\infty}$ an \emph{infinite dimensional vector
bundle over }$J^{\infty}$ whose module of sections identifies naturally with
$\underline{\Omega}^{1}(J^{\infty},\pi_{\infty})$. We conclude that $\tau
_{0}^{\dag}:J^{\dag}\longrightarrow J^{\infty}$ is naturally interpreted as
the reduced multimomentum bundle of $\pi_{\infty}$. Denote by $\ldots
,x^{i},\ldots,u_{I}^{\alpha},\ldots,p_{\alpha}^{I.i},\ldots$ standard
coordinates on $J^{\dag}$.  We will also consider the bundle
structures $J_{k}^{\dag
}\longrightarrow M$, $k\geq0$, and $\tau^{\dag}:=\pi_{\infty}\circ\tau
_{0}^{\dag}:J^{\dag}\longrightarrow M$.

In the following we denote by $U^\prime \subset J^\infty$ a generic open subset. Notice that any (local) element $\vartheta
\in\mathscr{C}\Lambda^{1}\otimes\overline{\Lambda}{}^{n-1}=\underline{\Omega
}^{1}(J^{\infty},\pi_{\infty})$, in particular a (local) Legendre form, is
naturally interpreted as a section $\vartheta:U^{\prime}\longrightarrow
J^{\dag}$ of $\tau^\dag_0$. Put then
$\ldots,\vartheta_{\alpha}^{I.i}:=\vartheta^{\ast}(p_{\alpha}^{I.i}),\ldots$
which are local functions on $J^{\infty}$ such that $\vartheta=\vartheta
_{\alpha}^{I.i}(du_{I}^{\alpha}-u_{Ii}^{\alpha}dx^{i})\otimes d^{n-1}x_{i}$.
It follows that, locally,
\[
\overline{d}\vartheta=-(D_{i}\vartheta_{\alpha}^{I.i}+\delta_{Ji}^{I}%
\vartheta_{\alpha}^{J.i})(du_{I}^{\alpha}-u_{Ii}^{\alpha}dx^{i})\otimes
d^{n}x\in\mathscr{C}\Lambda^{1}\otimes\overline{\Lambda}{}^{n},
\]
where $\delta_{K}^{I}=0$ if $I\neq K$, while $\delta_{K}^{I}=1$ if $I=K$.

Now, in Example \ref{Example1}, put $\alpha=\pi_{\infty}:P=J^{\infty
}\longrightarrow M$ and $\nabla=\mathscr{C}$, the Cartan connection in
$\pi_{\infty}$. $\mathscr{L}\in\Lambda_{n}^{n}=\overline{\Lambda}{}^{n}$ is
then a Lagrangian density in $\pi$. Put $\Sigma:=\Sigma_{\mathscr{C}}$,
$\Theta_{\mathscr{L}}:=\Theta_{\mathscr{L},\mathscr{C}}$ and $\omega
_{\mathscr{L}}:=\omega_{\mathscr{L},\mathscr{C}}$. $\omega_{\mathscr{L}}$ is a
PD-Hamiltonian system on $\tau^{\dag}:J^{\dag}\longrightarrow M$ canonically
determined by $\mathscr{L}$. Locally,
\[
\omega_{\mathscr{L}}=dp_{\alpha}^{I.i}du_{I}^{\alpha}d^{n-1}x_{i}%
-dE_{\mathscr{L}}d^{n}x,
\]
where $E_{\mathscr{L}}:=p_{\alpha}^{I.i}u_{Ii}^{\alpha}-L$. Let $\sigma
:U\longrightarrow J^{\dag}$ be a local section of $\tau^{\dag}$, and $j:=\tau_{0}^{\dag}\circ\sigma:U\longrightarrow J^\infty$. Put
$\ldots,\sigma_{I}^{\alpha}:=\sigma^{\ast}(u_{I}^{\alpha})=j^{\ast}%
(u_{I}^{\alpha}),\ldots,\sigma_{\alpha}^{I.i}:=\sigma^{\ast}(p_{\alpha}%
^{I.i}),\ldots$ which are local functions on $M$. Then, locally,
\[
i_{j_1 \sigma}\omega_{\mathscr{L}}|_{\sigma}=[(-\sigma_{\alpha}^{I.i}%
,_{i}-\delta_{Ji}^{I}\sigma_{\alpha}^{J.i}+\partial_{\alpha}^{I}L\circ
j)d^{V}\!u_{\alpha}^{I}|_{\sigma}+(\sigma_{I}^{\alpha},_{i}-\sigma_{Ii}^{\alpha
})d^{V}\!p_{\alpha}^{I.i}|_{\sigma}]\otimes d^{n}x,
\]
and the PD-Hamilton equations determined by $\omega_{\mathscr{L}}$ read
locally
\[
\left\{
\begin{array}
[c]{l}%
p_{\alpha}^{I.i},_{i}=\partial_{\alpha}^{I}L-\delta_{Ji}^{I}\,p_{\alpha}%
^{J.i}\\
u_{I}^{\alpha},_{i}=u_{Ii}^{\alpha}%
\end{array}
\right.  .
\]
We call such equations the \emph{Euler-Lagrange-Hamilton (ELH) equations
determined by the Lagrangian density }$\mathscr{L}$. Notice that they are
first order PDEs (with an infinite number of dependent variables). Denote by
$\mathscr{E}_{ELH}\subset J^{\infty}\tau^{\dag}$ their infinite prolongation.
In the following theorem we characterize solutions of $\mathscr{E}_{ELH}$. As
a byproduct, we derive the relationship between the ELH equations and the EL equations.

\begin{thm}
\label{Theorem1}A local section $\sigma:U\longrightarrow J^{\dag}$ of
$\tau^{\dag}$ is a solution of the ELH equations
determined by the Lagrangian density $\mathscr{L}$ iff it is locally of the
form $\sigma=\vartheta\circ j^{\infty}s$ where 1) $s:U\longrightarrow E$ is a
solution of the EL equations $\mathscr{E}_{EL}$ and 2) $\vartheta:U^{\prime
}\longrightarrow J^{\dag}$ is a Legendre form for $\mathscr{L}$.
\end{thm}

\begin{pf}
Let $\sigma:U\longrightarrow J^{\dag}$ be a local section of $\tau^{\dag}$. First of all, let $\sigma$ be of the form
$\sigma=\vartheta\circ j$ where 1) $j:U\longrightarrow J^{\infty}$ is a local
section of $\pi_{\infty}$ and 2) $\vartheta:U^{\prime}\longrightarrow J^{\dag
}$ is a local section of $\tau_{0}^{\dag}:J^{\dag}\longrightarrow J^{\infty}$. Then,
\[
\sigma_{\alpha}^{I.i},_{i}=D_{i}\vartheta_{\alpha}^{I.i}\circ j.
\]
Therefore, locally,
\begin{align*}
i_{j_1 \sigma}\omega_{\mathscr{L}}|_{\sigma}  &  =[[(-D_{i}\vartheta
_{\alpha}^{I.i}-\delta_{Ji}^{I}\vartheta_{\alpha}^{J.i}+\partial_{\alpha}%
^{I}L)\circ j]d^{V}\!u_{\alpha}^{I}|_{j}+(j_{I}^{\alpha},_{i}-j_{Ii}^{\alpha
})d^{V}\!p_{\alpha}^{I.i}|_{\sigma}]\otimes d^{n}x\\
&  =(\overline{d}\vartheta+d^{V}\!\mathscr{L})|_{j}+(j_{I}^{\alpha},_{i}%
-j_{Ii}^{\alpha})d^{V}\!p_{\alpha}^{I.i}|_{\sigma}\otimes d^{n}x,
\end{align*}
where $\ldots,j_{I}^{\alpha}:=j^{\ast}(u_{I}^{\alpha}),\ldots$ and they are
local functions on $M$. Thus, if $\vartheta$ is a Legendre form and
$j=j_{\infty}s$ for some local solution $s:U\longrightarrow E$ of the EL
equations then, in particular, $j_{I}^{\alpha}{},_{i}=j_{Ii}^{\alpha}$, $\alpha=1,\ldots,m$, $i=1,\ldots,n$, and
\[
i_{j_1 \sigma}\omega_{\mathscr{L}}|_{\sigma}=(\overline{d}\vartheta
+d^{V}\!\mathscr{L})|_{j}+(j_{I}^{\alpha},_{i}-j_{Ii}^{\alpha})d^{V}\!p_{\alpha
}^{I.i}|_{\sigma}\otimes d^{n}x=\boldsymbol{E}(\mathscr{L})|_{j}=0.
\]
On the other hand, let $\sigma:U\longrightarrow J^{\dag}$ be a local section
of $\tau^{\dag}$ and $j:=\tau_{0}^{\dag}\circ\sigma:U\longrightarrow
J^{\infty}$. Locally, there always exists a section $\vartheta:U^{\prime
}\longrightarrow J^{\dag}$ of $\tau_{0}^{\dag}$, such that $\sigma
=\vartheta\circ j$. Notice, preliminarily, that $\vartheta$ is not uniquely
determined by $\sigma$ except for its restriction to $\operatorname{im}j$. If
$\sigma$ is a solution of the ELH equations then, locally,
\[
0=i_{j_1 \sigma}\omega_{\mathscr{L}}|_{\sigma}=(\overline{d}\vartheta
+d^{V}\!\mathscr{L})|_{j}+(j_{I}^{\alpha},_{i}-j_{Ii}^{\alpha})d^{V}\!p_{\alpha
}^{I.i}|_{\sigma}\otimes d^{n}x.
\]
Since $(d^{V}\!p_{\alpha}^{I.i})|_{\sigma}\otimes d^{n}x$ and $(\overline
{d}\vartheta+d^{V}\!\mathscr{L})|_{j}$ are linearly independent, it follows
that
\[
\left\{
\begin{array}
[c]{l}%
(\overline{d}\vartheta+d^{V}\!\mathscr{L})|_{j}=0\\
j_{I}^{\alpha},_{i}=j_{Ii}^{\alpha}.
\end{array}
\right.  .
\]
In particular, $j=j_{\infty}s$, where $s=\pi_{\infty,0}\circ j$.

Now, let $\vartheta_{0}$ be a Legendre form for $\mathscr{L}$. Then
$d^{V}\!\mathscr{L}=\boldsymbol{E}(\mathscr{L})-\overline{d}\vartheta_{0}$ and,
therefore, $(\overline{d}\vartheta-\overline{d}\vartheta_{0}+\boldsymbol{E}%
(\mathscr{L}))|_{j}=0$. Recall that $\overline{d}$ restricts to $j=j_{\infty
}s$ (Remark \ref{Remark1}). Thus,
\[
\overline{d}|_{j}(\vartheta-\vartheta_{0})|_{j}=\boldsymbol{E}%
(\mathscr{L})|_{j}.
\]
In particular, $\boldsymbol{E}(\mathscr{L})|_{j}$ is $\overline{d}|_{j}%
$-exact. In view of Remark \ref{Remark1}, this is only possible if
$\boldsymbol{E}(\mathscr{L})|_{j}=0$, i.e., $s$ is a solution of the EL
equations. We conclude that
\[
\overline{d}|_{j}(\vartheta-\vartheta_{0})|_{j}=0,
\]
i.e., $(\vartheta-\vartheta_{0})|_{j}$ is $\overline{d}|_{j}$-closed. Again in
view of Remark \ref{Remark1}, this shows that, locally,
\[
(\vartheta-\vartheta_{0})|_{j}=\overline{d}|_{j}\nu|_{j}=\overline{d}\nu|_{j}%
\]
for some $\nu\in\mathscr{C}\Lambda^{1}\otimes\overline{\Lambda}{}^{n-2}$. In
particular, we can put $\vartheta=\vartheta_{0}+\overline{d}\nu$ and,
therefore, $\vartheta$ is a Legendre form for $\mathscr{L}$ as well.
\end{pf}

We now prove a \emph{formal version} of the above theorem. Put $p:=\tau
_{\infty,0}^{\dag}\circ\tau_{0}^{\dag}:J^{\infty}\tau^{\dag}\longrightarrow
J^{\infty}$.

\begin{thm}
\label{Theorem2}$p(\mathscr{E}_{ELH})\subset\mathscr{E}_{EL}$ and
$p:\mathscr{E}_{ELH}\longrightarrow\mathscr{E}_{EL}$ is a covering of PDEs.
\end{thm}

\begin{pf}
In $J^{\infty}\tau^{\dag}$ consider the submanifold $\mathscr{E}_{L}$ made of
$\infty$th jets of (local) sections $\sigma:U\longrightarrow J^{\dag}$ of the form $\sigma=\vartheta\circ j_{\infty}s$,
where $s:U\longrightarrow E$ is a local section of $\pi$, and $\vartheta
:U^{\prime}\longrightarrow J^{\dag}$ is a local Legendre form. It can be easily checked that
$\mathscr{E}_{L}$ is locally defined by
\begin{equation}
\left\{
\begin{array}
[c]{l}%
p_{\alpha}^{I.i}{}_{|Ki}+\delta_{Ji}^{I}\,p_{\alpha}^{J.i}{}{}_{|K}%
=D_{K}(\partial_{\alpha}^{I}L)-\delta_{\mathsf{O}}^{I}D_{K}\tfrac{\delta
L}{\delta u^{\alpha}}\\
u_{I}^{\alpha}{}_{|K}=u_{IK}^{\alpha}%
\end{array}
\right.  . \label{Eq1}%
\end{equation}
Clearly, the Cartan distribution restricts to $\mathscr{E}_{L}$ and,
therefore, $\mathscr{E}_{L}$ can be interpreted as a PDE. Moreover, it is
easily seen from (\ref{Eq1}) that $\mathscr{E}_{L}$ covers $J^{\infty}$ via
$p$. Denote by
\[
D_{j}^{\prime}=\tfrac{\partial}{\partial x^{j}}+u_{I}^{\alpha}{}_{|Jj}\tfrac{\partial}{\partial
u_{I}^{\alpha}{}_{|J}}+p_{\alpha}^{I.i}{}_{|Jj}\tfrac{\partial}{\partial
p_{\alpha}^{I.i}{}_{|J}}%
\]
the $j$th total derivative on $J^{\infty}\tau^{\dag}$, $j=1,\ldots,n$.
$\mathscr{E}_{ELH}$ is locally defined by
\begin{equation}
\left\{
\begin{array}
[c]{l}%
p_{\alpha}^{I.i}{}_{|Ki}=D_{K}^{\prime}(\partial_{\alpha}^{I}L)-\delta
_{Ji}^{I}\,p_{\alpha}^{J.i}{}{}_{|K}\\
u_{I}^{\alpha}{}_{|Ki}=u_{Ii}^{\alpha}{}_{|K}%
\end{array}
\right.  , \label{Eq}%
\end{equation}
which is equivalent to%
\[
\left\{
\begin{array}
[c]{l}%
p_{\alpha}^{I.i}{}_{|Ki}=D_{K}(\partial_{\alpha}^{I}L)-\delta_{Ji}%
^{I}\,p_{\alpha}^{J.i}{}_{|K}\\
u_{I}^{\alpha}{}_{|K}=u_{IK}^{\alpha}{}%
\end{array}
\right.  .
\]
Moreover, on $\mathscr{E}_{ELH}$
\[
(-)^{|I|}p_{\alpha}^{I.i}{}_{|KIi}=D_{K}\tfrac{\delta L}{\delta u^{\alpha}%
}-(-)^{|I|}\delta_{Ji}^{I}\,p_{\alpha}^{J.i}{}_{|KI}=D_{K}\tfrac{\delta
L}{\delta u^{\alpha}}+(-)^{|I|}p_{\alpha}^{I.i}{}_{|KIi},
\]
and, therefore, $D_{K}\tfrac{\delta L}{\delta u^{\alpha}}=0$, $\alpha=1,\ldots,m$. It then follows from (\ref{Eq1}), that
$\mathscr{E}_{ELH}=\mathscr{E}_{L}\cap p^{-1}(\mathscr{E}_{EL})$. In
particular, $p:\mathscr{E}_{ELH}\longrightarrow\mathscr{E}_{EL}$ is the
\textquotedblleft restriction\textquotedblright\ of $p:\mathscr{E}_{L}%
\longrightarrow J^{\infty}$ to $\mathscr{E}_{EL}\subset J^{\infty}$ and,
therefore, is a covering.
\end{pf}

\section{Natural Transformations of the Euler-Lagrange-Hamilton
Equations\label{SecTransfELH}}

Properties of Legendre forms discussed in Remark \ref{Remark2} correspond to
specific properties of the ELH equations which we discuss in this section.

First of all, notice that the ELH equations are canonically associated to a
Lagrangian density. But, how do the ELH equations change when changing the
Lagrangian density into a $\overline{d}$-cohomology class? In particular, does
an action functional uniquely determine a system of ELH equations or not? In
order to answer these questions consider $\vartheta\in\mathscr{C}\Lambda
^{1}\otimes\overline{\Lambda}{}^{n-1}$. $\vartheta$ determines an automorphism
$\Psi_{\vartheta}:J^{\dag}\longrightarrow J^{\dag}$ of the fiber bundle
$\tau_{0}^{\dag}$ via
\[
\Psi_{\vartheta}(P):=P-\vartheta_{\theta},\quad P\in J^{\dag},\,\theta
=\tau^{\dag}_0(P)\in J^{\infty}.
\]
In particular, $\tau_{0}^{\dag}\circ\Psi_{\vartheta}=\tau_{0}^{\dag}$.
Clearly, $\Psi_{\vartheta}^{-1}=\Psi_{-\vartheta}$.

\begin{lem}
\label{Lemma1}$\Psi_{\vartheta}^{\ast}(\omega_{\mathscr{L}})=\omega
_{\mathscr{L}}-\tau_{0}^{\dag}{}^{\ast}(d\vartheta)$.
\end{lem}

\begin{pf}
Compute,
\begin{align*}
\Psi_{\vartheta}^{\ast}(\omega_{\mathscr{L}})  &  =\Psi_{\vartheta}^{\ast
}(d\Theta_{\mathscr{L}})\\
&  =d\Psi_{\vartheta}^{\ast}(\Theta_{\mathscr{L}})\\
&  =d[(\Psi_{\vartheta}^{\ast}\circ\Sigma^{\ast})(\Theta)+(\Psi_{\vartheta
}^{\ast}\circ\tau_{0}^{\dag}{}^{\ast})(\mathscr{L})]\\
&  =d[(\Psi_{\vartheta}^{\ast}\circ\Sigma^{\ast})(\Theta)+(\tau_{0}^{\dag
}\circ\Psi_{\vartheta})^{\ast}(\mathscr{L})]\\
&  =d[(\Psi_{\vartheta}^{\ast}\circ\Sigma^{\ast})(\Theta)+\tau_{0}^{\dag}%
{}^{\ast}(\mathscr{L})].
\end{align*}
Now, since, locally, $\ldots,\Psi_{\vartheta}^{\ast}(p_{\alpha}^{I.i}%
)=p_{\alpha}^{I.i}-\vartheta_{\alpha}^{I.i},\ldots$, we have
\begin{align*}
(\Psi_{\vartheta}^{\ast}\circ\Sigma^{\ast})(p_{\alpha}^{I.i})  &  =p_{\alpha
}^{I.i}-\vartheta_{\alpha}^{I.i},\\
(\Psi_{\vartheta}^{\ast}\circ\Sigma^{\ast})(p)  &  =(p_{\alpha}^{I.i}%
-\vartheta_{\alpha}^{I.i})u_{Ii}^{\alpha}.
\end{align*}
Thus, locally
\begin{align*}
(\Psi_{\vartheta}^{\ast}\circ\Sigma^{\ast})(\Theta)  &  =(p_{\alpha}%
^{I.i}-\vartheta_{\alpha}^{I.i})du_{I}^{\alpha}d^{n-1}x_{i}-(p_{\alpha}%
^{I.i}-\vartheta_{\alpha}^{I.i})u_{Ii}^{\alpha}d^{n}x\\
&  =\Sigma^{\ast}(\Theta)-\tau_{0}^{\dag}{}^{\ast}(\vartheta).
\end{align*}
We conclude that
\begin{align*}
\Psi_{\vartheta}^{\ast}(\omega_{\mathscr{L}})  &  =d[(\Psi_{\vartheta}^{\ast
}\circ\Sigma^{\ast})(\Theta)+\tau_{0}^{\dag}{}^{\ast}(\mathscr{L})]\\
&  =d[\Sigma^{\ast}(\Theta)-\tau_{0}^{\dag}{}^{\ast}(\vartheta)+\tau_{0}%
^{\dag}{}^{\ast}(\mathscr{L})]\\
&  =\omega_{\mathscr{L}}-\tau_{0}^{\dag}{}^{\ast}(d\vartheta).
\end{align*}

\end{pf}

\begin{thm}
\label{TheorIso}Let $\mathscr{L}^{\prime}=\mathscr{L}+\overline{d}\varrho$,
$\varrho\in\overline{\Lambda}{}^{n-1}$, be another Lagrangian density (thus,
$\mathscr{L}^{\prime}$ determines the same EL equations as $\mathscr{L}$).
Then $\Psi_{d^{V}\!\varrho}^{\ast}(\omega_{\mathscr{L}})=\omega
_{\mathscr{L}^{\prime}}$.
\end{thm}

\begin{pf}
Notice, preliminarily, that
\begin{align*}
\tau_{0}^{\dag}{}^{\ast}(dd^{V}\!\varrho)  &  =\tau_{0}^{\dag}{}^{\ast
}(\overline{d}d^{V}\!\varrho)\\
&  =-\tau_{0}^{\dag}{}^{\ast}(d^{V}\!\overline{d}\varrho)\\
&  =-\tau_{0}^{\dag}{}^{\ast}(d\overline{d}\varrho)\\
&  =-d\tau_{0}^{\dag}{}^{\ast}(\overline{d}\varrho).
\end{align*}
Therefore, in view of the above lemma,
\begin{align*}
\Psi_{d^{V}\!\varrho}^{\ast}(\omega_{\mathscr{L}})  &  =\omega_{\mathscr{L}}%
-\tau_{0}^{\dag}{}^{\ast}(dd^{V}\!\varrho)\\
&  =d[\Sigma^{\ast}(\Theta)+\tau_{0}^{\dag}{}^{\ast}(\mathscr{L})]+d\tau
_{0}^{\dag}{}^{\ast}(\overline{d}\varrho)\\
&  =d[\Sigma^{\ast}(\Theta)+\tau_{0}^{\dag}{}^{\ast}(\mathscr{L}+\overline
{d}\varrho)]\\
&  =d\Theta_{\mathscr{L}^{\prime}}\\
&  =\omega_{\mathscr{L}^{\prime}}.
\end{align*}

\end{pf}

\begin{cor}
An action $[\mathscr{L}]\in\overline{H}{}^{n}$, $\mathscr{L}\in\overline
{\Lambda}{}^{n}$, uniquely determines a system of ELH equations, modulo
isomorphisms of PD-Hamiltonian systems.
\end{cor}

We conclude that the ELH equations are basically determined by the sole action
functional and not a specific Lagrangian density.

\begin{thm}
Let $\vartheta\in\mathscr{C}\Lambda^{1}\otimes\overline{\Lambda}{}^{n-1}$ be
$\overline{d}$-closed, hence $\overline{d}$-exact. Then, for every Lagrangian
density $\mathscr{L}\in\overline{\Lambda}{}^{n}$, $\Psi_{\vartheta}$ is a
symmetry of the ELH equations determined by $\mathscr{L}$ in the sense that
$j_{\infty}\Psi_{\vartheta}:J^{\infty}\tau^{\dag}\longrightarrow J^{\infty
}\tau^{\dag}$ preserves $\mathscr{E}_{ELH}$.
\end{thm}

\begin{pf}
By definition of infinite prolongations of a PDE and infinite prolongation of
a morphism of bundles, it is enough to prove that $j_{1}\Psi_{\vartheta}%
:J^{1}\tau^{\dag}\longrightarrow J^{1}\tau^{\dag}$ preserves
$\mathscr{E}_{ELH}^{(0)}:=\mathscr{E}_{\omega_{\mathscr{L}}}^{(0)}\subset
J^{1}\tau^{\dag}$. Notice, preliminarily, that, in view of the proof of
Theorem \ref{Theorem1}, we have
\[
(j_{1}\tau_{0}^{\dag})(\mathscr{E}_{ELH}^{(0)})\subset\operatorname{im}%
\mathscr{C}\subset J^{1}\pi_{\infty}.
\]
Now, let $c\in\mathscr{E}_{ELH}^{(0)}$, $P:=\tau_{1,0}^{\dag}(c)$ and $\xi\in
T_{P}J^{\dag}$ be a tangent vector, vertical with respect to $\tau^{\dag}$.
Consider also $c^{\prime}:=(j_{1}\Psi_{\vartheta})(c)$, $P^{\prime}%
:=\Psi_{\vartheta}(P)=\tau_{1,0}^{\dag}(c^{\prime})$ and $\xi^{\prime}%
:=d\Psi_{\vartheta}(\xi)$. In particular, $\xi^{\prime}\in T_{P^{\prime}%
}J^{\dag}$ is vertical with respect to $\tau^{\dag}$ as well. Let us prove
that $c^{\prime}\in\mathscr{E}_{ELH}^{(0)}$. In view of Lemma \ref{Lemma1},
\[
\Psi_{\vartheta}^{\ast}(\omega_{\mathscr{L}})=\omega_{\mathscr{L}}-\tau
_{0}^{\dag}{}^{\ast}(d\vartheta)=\omega_{\mathscr{L}}-\tau_{0}^{\dag}{}^{\ast
}(d^{V}\!\vartheta).
\]
Compute
\[
i_{\xi^{\prime}}i_{c^{\prime}}(\omega_{\mathscr{L}})_{P^{\prime}}%
=i_{\xi}i_{c}\Psi_{\vartheta}^{\ast}(\omega_{\mathscr{L}}%
)_{P}=i_{\xi}i_{c}(\omega_{\mathscr{L}})_{P}-i_{\xi}i_{c} [\tau^\dag_0{}^\ast(d^V\!\vartheta)]_P=-i_{\xi^{\prime\prime}}%
i_{\mathscr{C}_{\theta}}(d^{V}\!\vartheta)_{\theta}=0,
\]
where $\theta=\tau_{0}^{\dag}(P)\in J^{\infty}$ and $\xi^{\prime\prime}%
=(d\tau_{0}^{\dag})(\xi)\in T_{\theta}J^{\infty}$ is a tangent vector, vertical
with respect to $\pi_{\infty}$. It follows from the arbitrariness of
$\xi^{\prime}$, that $i_{c^{\prime}}(\omega_{\mathscr{L}})_{P^{\prime}}=0$.
\end{pf}

\section{Hamiltonian Formalism\label{SecHamForm}}

In this section we present our proposal of an Hamiltonian formalism for higher
order Lagrangian field theories. Such proposal is free from ambiguities in
that it depends only on the choice of a Lagrangian density and its order.
Moreover, $\overline{d}$-cohomologous Lagrangians of the same order determine equivalent
\textquotedblleft Hamiltonian theories\textquotedblright.

First of all, we define a \textquotedblleft finite dimensional
version\textquotedblright\ of the ELH equations (see also \cite{c...09}). In order to do this, notice
that, in view of Remark \ref{Remark3}, for all $k\geq 0$, $W_{k}$ is canonically isomorphic to
the $C^{\infty}(J^{k+1})$-module of sections of the induced bundle
$\pi_{k+1,k}^{\circ}(J^{\dag}\pi
_{k})\longrightarrow J^{k+1}$. We conclude that $J_{k}^{\dag}\longrightarrow J^{\infty}$ is canonically isomorphic to the
pull-back bundle $\pi_{\infty,k}^{\circ}(J^{\dag}\pi_{k})\longrightarrow J^{\infty}$, $k\geq0$. Notice that the coordinates $\ldots,p_{\alpha}^{I.i},\ldots$,
$|I|{}\leq k$, on $J_k^\dag$ identify with the pull-backs of the corresponding natural coordinates on $J^{\dag}\pi_{k}$ which we again denote by $\ldots,p_{\alpha}^{I.i},\ldots$.

Now, let $\mathscr{L}\in\overline{\Lambda}{}^{n}$ be a Lagrangian density of
order $l+1$, i.e., $\mathscr{L}\in\overline{\Lambda}{}^{n}\cap
\Lambda(J^{l+1})$. Let $\omega_{l}^{\prime}$ be the pull-back of $\omega
_{\mathscr{L}}$ onto $J^\dag_l$. $\omega
_{l}^{\prime}$ is a PD-Hamiltonian system on $J^\dag_l \longrightarrow M$, and it is
locally given by
\[
\omega_{l}^{\prime}=\sum_{|I|{}\leq l}dp_{\alpha}^{I.i}du_{I}^{\alpha}%
d^{n-1}x_{i}-dE_{l}d^{n}x,
\]
where $E_{l}=\sum_{|I|{}\leq l}p_{\alpha}%
^{I.i}u_{Ii}^{\alpha}-L$ is the restriction of $E_\mathscr{L}$ to $J^\dag_l$. Notice that $\omega_{l}^{\prime}$ is also the pull-back via $J^\dag_l \longrightarrow \pi^\circ_{l+1,l}(J^\dag \pi_l)$ of a (unique) PD-Hamiltonian system $\omega_{l}$ on $\pi^\circ_{l+1,l}(J^\dag \pi_l) \longrightarrow M$. $\omega_l$ is locally given by the same formula as $\omega_{l}^{\prime}
$ and it is a constrained PD-Hamiltonian system, i.e., its first
constraint bundle $\mathscr{P}\longrightarrow M$ is a proper subbundle of
$\pi^\circ_{l+1,l}(J^\dag \pi_l) \longrightarrow M$. Let us compute it. Let $P\in\pi^\circ_{l+1,l}(J^\dag \pi_l)$ and
$\theta:=\pi^\circ_{l+1,l}(\tau_0^\dag \pi_l)(P)\in J^{l+1}$. Then $P\in\mathscr{P}$
iff there exists $c$ in the first jet bundle of $\pi^\circ_{l+1,l}(J^\dag \pi_l) \longrightarrow M$ such that $i_{c}(\omega_{l})_{P}=0$,
i.e., iff there exist real numbers $\ldots,c_{I}^{a}._{i},\ldots,c_{\alpha
}^{I.i}{}._{j},\ldots$, $|I|{}\leq l$, such that
\[
\left\{
\begin{array}
[c]{ll}%
c_{\alpha}^{I.i}{}._{i}=(\partial_{\alpha}^{I}L)(\theta)-\delta_{Ji}%
^{I}\,P_{\alpha}^{J.i}{}, & |I|{}\leq l+1\\
c_{J}^{\alpha}._{i}=P_{Ji}^{\alpha}, & |J|{}\leq l
\end{array}
\right.
\]
where we put $c_{\alpha}^{I.i}._{i}=0$ for $|I|{}=l+1$, and $\ldots
,P_{Ji}^{\alpha}:=u_{Ji}^{\alpha}(P),\ldots,P^{K.i}_\alpha := p^{K.i}_\alpha (P),\ldots$, $|J|,|K|{}\leq l$, $\alpha
=1,\ldots,m$. Thus, for $|I|{}=l+1$, $P$ should be a solution of the system
\begin{equation}
\partial_{\alpha}^{I}L-\delta_{Ji}^{I}\,p_{\alpha}^{J.i}{}=0,\quad|I|{}=l+1.
\label{Eq3}%
\end{equation}
Equations (\ref{Eq3}) define $\mathscr{P}$ locally.
\begin{rem}\label{RemDim}
$\mathscr{P}$ is a submanifold of $\pi^\circ_{l+1,l}(J^\dag \pi_l)$ of the same dimension as $J^\dag \pi_l$, and $\mathscr{P} \longrightarrow J^{l+1}$ an affine subbundle of $\pi^\circ_{l+1,l}(J^\dag \pi_l) \longrightarrow J^{l+1}$.
\end{rem}
Let $\mathscr{P}_0$ be the image of $\mathscr{P}$ under the projection $\pi^\circ_{l+1,l}(J^\dag \pi_l) \longrightarrow J^\dag \pi_{l}$.

\begin{assum}\label{AssReg}
\label{Assumption1}We assume $\mathscr{P}_{0}$ to be a submanifold of
$J^{\dag}\pi_{l}$ and $\tau^{\dag}\pi_{l}|_{\mathscr{P}_{0}}:\mathscr{P}_{0}%
\longrightarrow M$ to be a smooth subbundle of $\tau^{\dag}\pi_{l}$. We also assume that the projection $q:\mathscr{P}\longrightarrow\mathscr{P}_{0}$
is a smooth submersion with connected fibers.
\end{assum}

Notice that, as usual, all the above regularity conditions are true if we
restrict all the involved maps to suitable open subsets.

The following commutative diagram summarizes the above described picture:
\[%
\begin{array}
[c]{c}%
\xymatrix@C=30pt@R=35pt{ 
          J^\dag 
                      \ar[d]       & \   J^\dag_{l} 
          \ar@{_{(}->}[l] 
          \ar[dl] 
          \ar[d]    
                      &    \\  
          J^\infty 
                      \ar[d]       
                      & \   \pi^\circ_{l+1,l}(J^\dag\pi_l)
          \ar[dl]
          \ar[d]    
                      &     \ \mathscr{P}  
          \ar@{_{(}->}[l]
          \ar[d]^-{q}        \\
          J^{l+1} 
          \ar[d]
          & J^\dag \pi_{l} 
          \ar[ddl]^(.3){\,\tau^\dag \pi_{l}} 
          \ar[dl]_-{\tau_0^\dag \pi_{l}}  
          & \ \mathscr{P}_0  
          \ar@{_{(}->}[l]
          \ar[ddll] \\
          J^{l} 
          \ar[d]_-{\pi_{l}} 
          & 
          &  \\
          M & &  
}
\end{array}
.
\]

\begin{thm}
\label{Theorem3}Under the regularity Assumption \ref{Assumption1}, there
exists a unique PD-Hamiltonian system $\omega_{0}$ on $\mathscr{P}_{0}\longrightarrow M$, such that
$i_{\mathscr{P}}^{\ast}(\omega_{l})$ is the pull-back of $\omega_0$ via $q: \mathscr{P} \longrightarrow \mathscr{P}_0$.
\end{thm}

\begin{pf}
Since $q:\mathscr{P}\longrightarrow\mathscr{P}_{0}$ has
connected fibers and $i^\ast_\mathscr{P}(\omega)$ is a closed form, it is enough to prove that $i_{\overline{Y}}i_{\mathscr{P}}%
^{\ast}(\omega_{l})=0$ for
all vector fields $\overline{Y}\in\mathrm{D}(\mathscr{P})$ vertical with
respect to $q$. Let $Y\in\mathrm{D}(\pi^\circ_{l+1,l}(J^\dag\pi_l))$ be vertical with respect to
$\pi^\circ_{l+1,l}(J^\dag\pi_l) \longrightarrow J^\dag \pi_l$, and $\overline{Y}:=Y|_{\mathscr{P}}$. Then $\overline{Y}$ is locally of the form
\[
\overline{Y}=\sum_{|K|{}=l+1}Y_{K}^{\beta}\partial_{\beta}^{K}|_{\mathscr{P}}%
,
\]
for some $\ldots,Y_{K}^{\beta},\ldots$ local functions on $\mathscr{P}$. Now
$\overline{Y}$ $\in\mathrm{D}(\mathscr{P})$ iff, locally,
\[
\sum_{|I|{}=l+1}Y_{K}^{\beta}\partial_{\beta}^{K}\partial_{\alpha}%
^{I}L|_{\mathscr{P}}=0.
\]
Compute
\[
\overline{Y}(E_{l}|_{\mathscr{P}})=\sum_{|K|{}=l+1}Y_{K}^{\beta}%
\partial_{\beta}^{K}E_{l}|_{\mathscr{P}}=\sum_{|I|{}=l+1}Y_{I}^{\alpha}%
(\delta_{Ji}^{I}\,p_{\alpha}^{J}{}^{.i}-\partial_{\alpha}^{I}L)|_{\mathscr{P}}%
=0.
\]
Thus $E_l |_\mathscr{P}$ is the pull-back via $q$ of a (unique) local function $H$ on $\mathscr{P}_0$.
Moreover,
\[
i_{\overline{Y}}i_{\mathscr{P}}^{\ast}(\omega_{l})=-\overline{Y}%
(E_{l}|_{\mathscr{P}})d^{n}x=0.
\]
It follows from the arbitrariness of $\overline{Y}$ that $i_{\mathscr{P}}%
^{\ast}(\omega_{l})$ is the pull-back via $q$ of the PD-Hamiltonian system $\omega_0$ on $\mathscr{P}_0 \longrightarrow M$ locally defined as
\[
\omega_{0}=\sum_{|I|{}\leq l}i_{\mathscr{P}_{0}}^{\ast}(dp_{\alpha}%
^{I.i}du_{I}^{\alpha})d^{n-1}x_{i}-dHd^{n}x.
\]
\end{pf}

\begin{defn}
$\omega_{0}$ is called the \emph{PD-Hamiltonian system determined by the
}$(l+1)$\emph{th order Lagrangian density} $\mathscr{L}$, and the corresponding
PD-Hamilton equations are the \emph{Hamilton-de Donder-Weyl (HDW) equations
determined by }$\mathscr{L}$.
\end{defn}

\begin{defn}
A Lagrangian density $\mathscr{L}$ of order $l+1$ is \emph{regular at the
order} $l+1$ iff the map $\mathscr{P} \longrightarrow J^\dag \pi_l$ has maximum rank.
\end{defn}

The Lagrangian density $\mathscr{L}$ of order $l+1$ is regular at the order $l+1$
iff the matrix%
\[
\mathbf{H}(L)(\theta):=\left\Vert (\partial_{\beta}^{K}\partial_{\alpha}^{I}L)(\theta)\right\Vert
{}_{(\beta,K)}^{(\alpha,I)},\quad|I|,|K|{}=l+1,
\]
where the pairs $(\alpha,I)$ and $(\beta,K)$ are understood as single indexes,
has maximum rank at every point $\theta\in J^{l+1}$. In its turn, this implies that $\mathscr{P}_{0}$ is an open submanifold of $J^{\dag}\pi_{l}$ and, in view of Remark \ref{RemDim} and Assumption \ref{AssReg}, $q: \mathscr{P} \longrightarrow \mathscr{P}_0$ is a diffeomorphism. In particular, $\omega_{0}$ is a PD-Hamiltonian system on an open subbundle of $\tau^{\dag
}\pi_{l}$ locally given by
\[
\omega_{0}=\sum_{|I|{}\leq l}dp_{\alpha}^{I.i}du_{I}^{\alpha}d^{n-1}%
x_{i}-dHd^{n}x,
\]
where, now, $H$ is a local function on $J^{\dag}\pi_{l}$. In this case, as expected,
the HDW equations read locally
\[
\left\{
\begin{array}
[c]{ll}%
p_{\alpha}^{I.i},_{i}=-\tfrac{\partial H}{\partial u_{I}^{\alpha}}\\
u_{I}^{\alpha},_{i}=\tfrac{\partial H}{\partial p_{\alpha}^{I.i}}%
\end{array}
\right.  .
\]
Notice that the HDW equations are canonically associated to a Lagrangian
density and its order and no additional structure is required to define them.
Moreover, in view of Theorem \ref{TheorIso}, two Lagrangian densities of the
same order determining the same system of EL equations, also determine equivalent
HDW equations. Finally, to write down the HDW equations there is no need of
a distinguished Legendre transform. Actually, the emergence of ambiguities
in all Hamiltonian formalisms for higher order field theories proposed in the
literature seems to rely on the common attempt to define first a higher order
analogue of the Legendre transform and, only thereafter, the
\textquotedblleft Hamiltonian theory\textquotedblright. In the next section we
present our own point of view on the Legendre transform in higher order
Lagrangian field theories.

\section{\label{SecLegTransf}The Legendre Transform}

Keeping the same notations as in the previous section, denote by ${}%
^{l}\!\mathscr{E}_{ELH}$ the infinite prolongation of
the PD-Hamilton equations determined by  $\omega_{l}$ and by $p^{\prime}:\pi_{l+1,l}^\circ(J^\dag \pi_l)\longrightarrow E$ the natural projection.

\begin{prop}
\label{PropCov1}$(j_{\infty}p^{\prime})({}^{l}\!\mathscr{E}_{ELH}%
)\subset\mathscr{E}_{EL}$ and $j_{\infty}p^{\prime}:{}^{l}\!\mathscr{E}_{ELH}%
\longrightarrow\mathscr{E}_{EL}$ is a covering.
\end{prop}

\begin{pf}
The proof is the finite dimensional version of the proof of Theorem
\ref{Theorem2} and will be omitted (see also \cite{c...09}).
\end{pf}

Denote also by $\mathscr{E}_{H}^{\mathscr{P}}$ the infinite prolongation of the PD-Hamilton equations determined by 
$i_{\mathscr{P}}^{\ast}(\omega_{l})$ and by $\mathscr{E}_{H}$ the infinite prolongation of the HDW equations.

\begin{prop}
\label{PropCov2}$(j_{\infty}q)({}\mathscr{E}_{H}^{\mathscr{P}}%
)\subset\mathscr{E}_{H}$ and $j_{\infty}q:\mathscr{E}_{H}^{\mathscr{P}}%
\longrightarrow\mathscr{E}_{H}$ is a covering.
\end{prop}

\begin{pf}
It immediately follows from Theorem \ref{Theorem3} and Proposition \ref{Prop2}.
\end{pf}

Notice that, in view of Propositions \ref{Prop1}, \ref{PropCov1} and
\ref{PropCov2}, there is a diagram of morphisms of PDEs,%
\begin{equation}%
\begin{array}
[c]{c}%
\xymatrix@C=30pt@R=30pt{ {}^l \! \mathscr{E}_{ELH}\ \ar@{^{(}->}[r] \ar[d]_-{j_\infty p^\prime} &\ \mathscr{E}_H ^{\mathscr{P}}  \ar[d]^-{\ j_\infty q} \\
                    \mathscr{E}_{EL} & \mathscr{E}_H 
                     }
\end{array}
, \label{DiagNLLeg}%
\end{equation}
whose vertical arrows are coverings. Therefore, the inclusion $^{l}%
\!\mathscr{E}_{ELH}\subset\mathscr{E}_{H}^{\mathscr{P}}$ may be understood as
a non local morphism of $\mathscr{E}_{EL}$ into $\mathscr{E}_{H}$. We
interpret such morphism as Legendre transform according to the following

\begin{defn}
We call diagram (\ref{DiagNLLeg}) the \emph{Legendre transform determined by
the Lagrangian density} $\mathscr{L}$.
\end{defn}

Any Legendre form of order $l$, $\vartheta:J^{\infty}\longrightarrow
J_{l}^{\dag}$ $\longrightarrow\pi_{l+1,l}^\circ(J^\dag \pi_l)$, determines a
section $j_{\infty}\vartheta|_{\mathscr{E}_{EL}}:\mathscr{E}_{EL}%
\longrightarrow{}^{l}\!\mathscr{E}_{ELH}$ of the covering $j_{\infty}%
p^{\prime}:{}^{l}\!\mathscr{E}_{ELH}\longrightarrow\mathscr{E}_{EL}$ and,
therefore, via composition with $j_{\infty}q$, a concrete map
$\mathscr{E}_{EL}\longrightarrow\mathscr{E}_{H}$. Nevertheless, among these
maps, there is no distinguished one.

We now prove that, if $\mathscr{L}$ is regular at the order $l+1$, then
$\mathscr{E}_{H}$ itself covers $\mathscr{E}_{EL}$. This result should be
interpreted as the higher order analogue of the theorem stating the
equivalence of EL equations and HDW equations for first order theories with
regular Lagrangian (see, for instance, \cite{gs73}). Let us first prove the following

\begin{lem}
If $\mathscr{L}$ is regular at the order $l+1$, then $^{l}\!\mathscr{E}_{ELH}%
=\mathscr{E}_{H}^{\mathscr{P}}$.
\end{lem}

\begin{pf}
The proof is in local coordinates. Let $\sigma:U\longrightarrow
\pi_{l+1,l}^\circ(J^\dag \pi_l)$ be a local section of $\pi_{l+1,l}^\circ(J^\dag \pi_l) \longrightarrow M$. Suppose $\operatorname{im}\sigma\subset\mathscr{P}$. Then,
locally,
\[
\partial_{\alpha}^{I}L\circ\sigma-\delta_{Jj}^{I}\sigma_{\alpha}^{J.j}%
=0,\quad|I|{}=l+1.
\]
Now, $i_{j_1 \sigma}\omega_{l}|_{\sigma}$ is locally given by
\begin{equation*}
 i_{j_1 \sigma}\omega_{l}|_{\sigma}
 =\left[
{\textstyle\sum\nolimits_{|I|{}\leq l+1}}
(-\sigma_{\alpha}^{I.i},_{i}-\delta_{Jj}^{I}\sigma_{\alpha}^{J.j}%
+\partial_{\alpha}^{I}L\circ\sigma)d^V \!u_{I}^{\alpha}+%
{\textstyle\sum\nolimits_{|I|{}\leq l}}
(\sigma_{I}^{\alpha},_{i}-\sigma_{Ii}^{\alpha})d^V \!p_{\alpha}^{I.i}\right]
|_{\sigma}\otimes d^{n}x.
\end{equation*}
As already outlined, the annihilator of $\mathrm{D}(\mathscr{P})$ in
$\Lambda^{1}(\pi_{l+1,l}^\circ(J^\dag \pi_l))|_{\mathscr{P}}$ is locally spanned
by 1-forms
\[
\lambda_{\alpha}^{I}:=d(\partial_{\alpha}^{I}L-\delta_{Jj}^{I}p_{\alpha}%
^{J.j})|_{\mathscr{P}},\quad|I|{}=l+1.
\]
Therefore, $i_{j_1 \sigma}i_{\mathscr{P}}^{\ast}(\omega_{l})|_{\sigma}=0$
iff, locally,
\begin{equation}
i_{j_1 \sigma}\omega_{l}|_{\sigma}=\sum_{|I|{}=l+1}f_{I}^{a}\underline{\lambda}{}_{\alpha
}^{I}|_{\sigma} \otimes d^n x, \label{Eq2}%
\end{equation}
for some local functions $\ldots,f_{I}^{\alpha},\ldots$ on $\operatorname{im}%
\sigma$, where 
\[
\underline{\lambda}_{\alpha}^{I}:= d^V\! (\partial_{\alpha}^{I}L-\delta_{Jj}^{I}p_{\alpha}%
^{J.j})|_{\mathscr{P}}={\textstyle\sum\nolimits_{|K|{}\leq l+1}}
(\partial_{\beta}^{K}\partial_{\alpha}^{I}L\,d^V \! u_{K}^{\beta}-\delta_{Jj}%
^{I}d^V \! p_{\alpha}^{Jj})  |_{\mathscr{P}},\quad|I|{}=l+1.
\] Equations (\ref{Eq2}) read
\[%
\begin{array}
[c]{r}%
\sum_{|I|{}\leq l+1}(-\sigma_{\alpha}^{I.i},_{i}+\partial_{\alpha}^{I}%
L\circ\sigma-\delta_{Jj}^{I}\sigma_{\alpha}^{J.j}-\sum_{|K|{}=l+1}f_{K}%
^{\beta}\partial_{\beta}^{K}\partial_{\alpha}^{I}L\circ\sigma)d^V \! u_{I}^{\alpha
}|_{\sigma}\\
+\sum_{|I|{}<l}(\sigma_{I,i}^{\alpha}-\sigma_{Ii}^{\alpha})d^V \! p_{\alpha}%
^{I.i}|_{\sigma}+\sum_{|I|{}=l}(\sigma_{I}^{\alpha},_{i}-\sigma_{Ii}^{\alpha
}+\tfrac{I[i]+1}{l+1}f_{Ii}^{\alpha})d^V \! p_{\alpha}^{I.i}|_{\sigma}=0
\end{array}
,
\]
where $I[i]$ is the number of times the index $i$ appears in the multiindex $I$.
Since the vertical forms $\dots,d^V \! u_{I}^{\alpha}|_{\sigma},\ldots,d^V \! p_{\alpha}%
^{I.i}|_{\sigma},\ldots$ are linearly independent, $i_{j_1 \sigma%
}i_{\mathscr{P}}^{\ast}(\omega_{l})|_{\sigma}=0$ iff, locally,
\begin{equation}
\left\{
\begin{array}
[c]{ll}%
-\sigma_{\alpha}^{I.i},_{i}+\partial_{\alpha}^{I}L\circ\sigma-\delta_{Jj}%
^{I}\sigma_{\alpha}^{J.j}-\sum_{|K|{}=l+1}f_{K}^{\beta}\partial_{\beta}%
^{K}\partial_{\alpha}^{I}L\circ\sigma=0, & |I|{}\leq l+1\\
\sigma_{I,i}^{\alpha}-\sigma_{Ii}^{\alpha}=0, & |I|{}<l\\
\sigma_{I}^{\alpha},_{i}-\sigma_{Ii}^{\alpha}+\tfrac{I[i]+1}{l+1}%
f_{Ii}^{\alpha}=0, & |I|{}=l
\end{array}
\right.  , \label{Eq4}%
\end{equation}
for some $\ldots,f_{I}^{\alpha},\ldots$. It follows from the third of
Equations (\ref{Eq4}) that
\begin{equation}
f_{Ii}^{\alpha}=-\tfrac{l+1}{I[i]+1}(\sigma_{I}^{\alpha},_{i}-\sigma
_{Ii}^{\alpha}),\quad|I|{}=l. \label{Eq6}%
\end{equation}
Moreover, since $\operatorname{im}\sigma\subset\mathscr{P}$, the first
equation, for $|I|{}=l+1$, gives
\[
0=\sum_{|K|{}=l+1}f_{K}^{\beta}\partial_{\beta}^{K}\partial_{\alpha}^{I}%
L\circ\sigma=\sum_{|J|{}=l}\tfrac{J[j]+1}{l+1}f_{Jj}^{\beta}\partial_{\beta
}^{Jj}\partial_{\alpha}^{I}L\circ\sigma=-\sum_{|J|{}=l}(\sigma_{J}^{\beta
},_{j}-\sigma_{Jj}^{\beta})\partial_{\beta}^{Jj}\partial_{\alpha}^{I}%
L\circ\sigma,
\]
and, in view of the regularity of $\mathscr{L}$ and Equations (\ref{Eq6}),
\[
\sigma_{I}^{\alpha},_{i}-\sigma_{Ii}^{\alpha}=f_{Ii}^{\alpha}=0,\quad|I|{}=l.
\]
Substituting again into (\ref{Eq4}), we finally find that the PD-Hamilton
equations $i_{j_1 \sigma}i_{\mathscr{P}}^{\ast}(\omega_{l})|_{\sigma}=0$ are
locally equivalent to equations
\[
\left\{
\begin{array}
[c]{ll}%
p_{\alpha}^{I.i},_{i}=\partial_{\alpha}^{I}L-\delta_{Jj}^{I}p_{\sigma}%
^{J.j}, & |I|{}\leq l+1\\
u_{I}^{\alpha},_{i}=u_{Ii}^{\alpha}, & |I|{}\leq l
\end{array}
\right.  ,
\]
which are the PD-Hamilton equations determined by $\omega_{l}$.
\end{pf}

Now, suppose that $\mathscr{L}$ is regular at
the order $l+1$. Then, as already mentioned in the previous section, $q: \mathscr{P} \longrightarrow \mathscr{P}_0$ is a diffeomorphism, and $q^\ast({\omega_0}) = i^\ast_\mathscr{P}(\omega_l)$. Therefore, $j_\infty q : \mathscr{E}^\mathscr{P}_H \longrightarrow \mathscr{E}_H$ is an isomorphism of PDEs and the Legendre transform (\ref{DiagNLLeg}) reduces to
\[%
\begin{array}
[c]{c}%
\xymatrix@C=30pt@R=30pt{ {}^l \! \mathscr{E}_{ELH}\ \ar@{=}[r] \ar[d]_-{j_\infty p^\prime} &\ \mathscr{E}_H ^{\mathscr{P}}  \ar[d]^-{\ j_\infty q} \ar@{}[d]^{\begin{sideways}$\widetilde{\quad\quad}$\end{sideways}} \\
                    \mathscr{E}_{EL} & \mathscr{E}_H 
                     }
\end{array}
.
\]
Moreover, $J^{\dag}\pi_{l}$ maps to $E$ via $\pi_{l,0}\circ\tau_{0}^{\dag}%
\pi_{l}$ and such map is a morphism of bundles (over $M$). The induced
morphism $J^{\infty}\tau_{0}^{\dag}\pi_{l}\longrightarrow J^{\infty
}$ restricts to a morphism of PDEs, $\kappa:\mathscr{E}_{H}%
\longrightarrow J^{\infty}$, locally defined as $\kappa^{\ast}(u_{K}^{\alpha
})=u_{\mathsf{O}}^{\alpha}{}_{|K}$, $|K|{}\geq0$. It is easy to show that
diagram
\[%
\begin{array}
[c]{c}%
\xymatrix@C=30pt@R=30pt{ {}^l \! \mathscr{E}_{ELH}\ \ar@{=}[rr] \ar[d]_-{j_\infty p^\prime}& &\ \mathscr{E}_H ^{\mathscr{P}}  \ar[d]^-{\ j_\infty q} \ar[d]^{\begin{sideways}$\widetilde{\quad\quad}$\end{sideways}} \\
                    \mathscr{E}_{EL}\ \ar@{^{(}->}[r] & J^\infty &\mathscr{E}_H \ar[l]_-{\kappa} 
                     }
\end{array}
,
\]
commutes, so that $\kappa = j_\infty p^\prime \circ (j_\infty q |_{\mathscr{E}_H})^{-1}$. Consequently, $\kappa(\mathscr{E}_{H})\subset\mathscr{E}_{EL}$ and
$\kappa:\mathscr{E}_{H}\longrightarrow\mathscr{E}_{EL}$ is a covering. Summarizing, we have
have proved the following

\begin{thm}
If $\mathscr{L}$ is regular at the order $l+1$, then $\mathscr{E}_{H}$ covers
$\mathscr{E}_{EL}$.
\end{thm}

Finally, it should be mentioned that in most cases, even if the Lagrangian density is not regular, $\mathscr{E}_{H}$ covers $\mathscr{E}_{EL}$ via $\kappa$ and, therefore, $\mathscr{E}_H ^{\mathscr{P}}$ itself covers $\mathscr{E}_{EL}$ (see the example in the next section). 

\section{An Example: The Korteweg-de Vries Action}\label{SecKdV}

The celebrated Korteweg-de Vries (KdV) equation
\begin{equation}
\phi_{t}-6\phi\phi_{x}+\phi_{xxx}=0\label{KdV}%
\end{equation}
can be derived from a variational principle as follows. Introduce the
\textquotedblleft potential\textquotedblright\ $u$ by putting $u_{x}=\phi$.
Equation (\ref{KdV}) becomes the fourth order non-linear equation
\begin{equation}
u_{tx}-6u_{x}u_{xx}+u_{xxxx}=0\label{KdVEL}%
\end{equation}
for sections of the trivial bundle $\pi:\mathbb{R}^{2}\times\mathbb{R}%
\ni(t,x;u)\longmapsto(t,x)\in\mathbb{R}^{2}$. In its turn, (\ref{KdVEL}) is
the EL equation determined by the action functional
\[
\int(u_{x}^{3}-\tfrac{1}{2}u_{x}u_{t}+\tfrac{1}{2}u_{xx}^{2})dtdx.
\]
Choose the second order Lagrangian density
\begin{equation}
\mathscr{L}=(u_{x}^{3}-\tfrac{1}{2}u_{x}u_{t}+\tfrac{1}{2}u_{xx}%
^{2})dtdx.\label{KdVLag}%
\end{equation}
Since the matrix
\[
\mathbf{H}(L) =\left(
\begin{array}
[c]{ccc}%
0 & 0 & 0\\
0 & 0 & 0\\
0 & 0 & 1
\end{array}
\right)
\]
has rank $1$, $\mathscr{L}$ is not regular. Let $t,x,u,u_{t},u_{x}%
,p^{.t},p^{.x},\ldots p^{i.j},\ldots$ be natural coordinates on
$J^{\dag}\pi_{1}$, $i,j=t,x$. Then
\[
\omega_{1}=dp^{.t}dudx-dp^{.x}dudt+dp^{t.t}du_{t}dx-dp^{t.x}du_{t}%
dt+dp^{x.t}du_{x}dx-dp^{x.x}du_{x}dt-dE_{\mathrm{KdV}}dtdx
\]
where
\[
E_{\mathrm{KdV}}:=p^{.t}u_{t}+p^{.x}u_{x}+p^{t.t}u_{tt}+(p^{t.x}%
+p^{x.t})u_{tx}+p^{x.x}u_{xx}-u_{x}^{3}+\tfrac{1}{2}u_{x}u_{t}-\tfrac{1}%
{2}u_{xx}^{2}.
\]
Accordingly, ${}^{1}\mathscr{E}_{ELH}$ reads
\begin{equation}
^{1}\mathscr{E}_{ELH}:\left\{
\begin{array}
[c]{ll}%
p^{.t},_{t}+p^{.x},_{x}=0 & \\
p^{t.t},_{t}+p^{t.x},_{x}=-\tfrac{1}{2}u_{x} & \\
p^{x.t},_{t}+p^{x.x},_{x}=3u_{x}^{2}-\tfrac{1}{2}u_{t} & \\
u,_{i}=u_{i} & i=t,x\\
u_{i},_{j}=u_{ij} & i,j=t,x\\
p^{t.t}=0 & \\
p^{t.x}+p^{x.t}=0 & \\
p^{x.x}-u_{xx}=0 &
\end{array}
\right.  .\label{KdVELH}%
\end{equation}
which clearly cover (\ref{KdVEL}). Notice that the last three equations in
(\ref{KdVELH}) define $\mathscr{P}$. Thus, $\mathscr{P}$ is coordinatized by
$t,x,u,u_{t},u_{x},u_{tt},u_{tx},p^{t.x},p^{x.x}$ and
\[
i_{\mathscr{P}}^{\ast}(\omega_{1})=dp^{.t}dudx-dp^{.x}dudt-dp^{t.x}%
(du_{t}dt+du_{x}dx)-dp^{x.x}du_{x}dt-dE_{\mathrm{KdV}}|_{\mathscr{P}}dtdx,
\]
where
\[
E_{\mathrm{KdV}}|_{\mathscr{P}}:=p^{.t}u_{t}+p^{.x}u_{x}+\tfrac{1}{2}%
(p^{x.x})^{2}-u_{x}^{3}+\tfrac{1}{2}u_{x}u_{t}.
\]
Accordingly, $\mathscr{E}_{H}^{\mathscr{P}}$ reads
\[
\mathscr{E}_{H}^{\mathscr{P}}:\left\{
\begin{array}
[c]{ll}%
p^{.t},_{t}+p^{.x},_{x}=0 & \\
p^{t.x},_{x}=-\tfrac{1}{2}u_{x} & \\
p^{t.x},_{t}+p^{x.x},_{x}=-3u_{x}^{2}+\tfrac{1}{2}u_{t} & \\
u,_{i}=u_{i} & i=t,x\\
u_{t},_{x}=u_{x},_{t} & \\
u_{x},_{x}=p^{x.x} &
\end{array}
\right.  .
\]
Notice that, even if the Lagrangian density is not regular, and variables $u_{tt},u_{tx}$
are undetermined, $\mathscr{E}_{H}^{\mathscr{P}}$ covers (\ref{KdVEL}). Finally, $\mathscr{P}$ is defined by the sixth and the
seventh equations in (\ref{KdVELH}) and, therefore, it is coordinatized by
$t,x,u,u_{t},u_{x},p^{t.x},p^{x.x}$. Thus, $\omega_{0}$ and $\mathscr{E}_{HDW}%
$ are given by exactly the same coordinate formulas as $i_{\mathscr{P}}^{\ast
}(\omega_{1})$ and $\mathscr{E}_{H}^{\mathscr{P}}$. In particular, $\mathscr{E}_{HDW}$ itself covers
$\mathscr{E}_{EL}$.

Finally, recall that the KdV equation is Hamiltonian, i.e., it can be
presented in the form $u_{t}=A(\boldsymbol{E}(\mathscr{H}))$, where
$\mathscr{H}$ is a top horizontal form in the infinite jet space of the bundle
$\mathbb{R}^{2}\ni(x;u)\longmapsto x\in\mathbb{R}$, and $A$ is a Hamiltonian
$\mathscr{C}$-differential operator (see, for instance, \cite{b...99}). Since
Hamiltonian PDEs play a prominent role in the theory of integrable systems, it
is worth to mention that such property (which is based on a 1+1,
\textquotedblleft covariance breaking\textquotedblright\ splitting of the
space of independent variables $(t,x)$) is directly related with the present
covariant Hamiltonian formalism as shown, for instance, in \cite{g88}. There
the author provides a multisymplectic framework for the KdV equation by
choosing, along the lines of \cite{k84b}, a \textquotedblleft
quasi-symmetric\textquotedblright\ Cartan form for the Lagrangian density
(\ref{KdVLag}). Such Cartan form is unique for a second order theory.
Therefore, the formalism of \cite{g88} is actually equivalent to ours, in the
special case of a second order theory.

\section*{Conclusions}

In this paper, using the geometric theory of PDEs, we solved the long standing
problem of finding a reasonably natural, higher order, field theoretic
analogue of Hamiltonian mechanics of Lagrangian systems. By naturality we mean
dependence on no structure other than the action functional. We achieved our
goal in two steps. First we found a higher order, field theoretic analogue of
the Skinner-Rusk mixed Lagrangian-Hamiltonian formalism \cite{s83,sr83,sr83b}
and, second, we showed that such theory projects naturally to a PD-Hamiltonian
system on a smaller space. The obtained Hamiltonian field equations enjoy the
following nice properties: 1) they are first order, 2) there is a canonical,
non-local embedding of the Euler-Lagrange equations into them, and 3) for
regular Lagrangian theories, they cover the Euler-Lagrange equations.
Moreover, for regular Lagrangian theories, the coordinate expressions of the
obtained field equations are nothing but the de Donder higher order field
equations. This proves that our theory is truly the coordinate-free
formulation of de Donder one \cite{d35}.

\end{document}